\documentclass{article}
\usepackage{amssymb}
\usepackage{amsfonts}
\usepackage{amsmath}

\begin{document}

\begin{center}
{\large Wajsberg algebras arising from binary block codes}

\bigskip

Cristina FLAUT and Radu VASILE%
\begin{equation*}
\end{equation*}
\end{center}

\textbf{Abstract.} {\small In this paper we presented some connections
between BCK-commutative bounded algebras, MV-algebras, Wajsberg algebras and
binary block codes. Using connections between these three algebras, we will
associate to each of them a binary block code and, in some circumstances, we
will prove that the converse is also true.}%
\begin{equation*}
\end{equation*}


\begin{equation*}
\end{equation*}

\bigskip 1. \textbf{Introduction}

\begin{equation*}
\end{equation*}%
\qquad\ \qquad

BCK-algebras were first introduced in mathematics by Y. Imai and K. Iseki,
in 1966, through the paper [II; 66], as a generalization of the concept of
set theoretic difference and propositional calculi. These algebras form an
important class of logical algebras and have many applications to various
domains of mathematics (group theory, functional analyses, sets theory,
etc.). Because of the necessity to establish certain rational logic systems
as a logical foundation for uncertain information processing, various types
of logical systems have been proposed. For this purpose, some logical
algebras appeared and have been researched.([WDH; 17]) One of these algebras
are MV-algebras, where MV is referred to \textquotedblright many
valued\textquotedblright ( [GA; 90]), which were originally introduced by
Chang in [CHA; 58]. He tried to provide a new proof for the completeness of
the \L ukasiewicz axioms for infinite valued propositional logic. These
algebras appeared in the specialty literature under equivalent names:
bounded commutative BCK-algebras or Wajsberg algebras, ([CT; 96]). Wajsberg
algebras were introduced in 1984, by Font, Rodriguez and Torrens, through
the paper \ [FRT; 84] as an alternative model for the infinite valued \L %
ukasiewicz propositional logic.

In the following, we present some connections between BCK-commutative
bounded algebras, MV-algebras, Wajsberg algebras and binary block codes and
we gave an algorithm to find all finite partial ordered Wajsberg algebras.%
{\small \ }These new approach allows us to find for these structures new and
interesting properties.%
\begin{equation*}
\end{equation*}%
\begin{equation*}
\end{equation*}

\begin{equation*}
\end{equation*}

\textbf{2. Preliminaries}%
\begin{equation*}
\end{equation*}

\textbf{Definition 2.1.} An algebra $(X,\ast ,\theta )$ of type $(2,0)$ is
called a \textit{BCI-algebra} if the following conditions are fulfilled:

$1)~((x\ast y)\ast (x\ast z))\ast (z\ast y)=\theta $, for all $x,y,z\in X;$

$2)~(x\ast (x\ast y))\ast y=\theta $, for all $x,y\in X;$

$3)~x\ast x=\theta $, for all $x\in X$;

$4)$ For all $x,y,z\in X$ such that $x\ast y=\theta ,y\ast x=\theta $, it
results $x=y$.

If a BCI-algebra $X$ satisfies the following identity:

$5)$ $\theta \ast x=\theta ,~$for all $x\in X,$ then $X$ is called a \textit{%
BCK-algebra}.

In a BCK-algebra we have the following order relation:%
\begin{equation*}
x\leq y\text{ if and only if \ }x\ast y=\theta .
\end{equation*}

If the algebra $\left( X,\ast ,\theta \right) $ has an element $1$ such that 
$x\ast 1=\theta ,$ for all $x\in X~$(that means $x\leq 1,$ for all $x\in X)$%
, then the BCK-algebra $X$ is called \textit{bounded}. If $x\wedge y=y\wedge
x,$ for all $x,y\in X$, where $x\wedge y=y\ast \left( y\ast x\right) $, for
all $x,y\in X$, then $X$ is called a \textit{commutative} BCK-algebra. For
other details regarding BCK algebras, the readers are referred to [AAT; 96],
[Me-Ju; 94].\medskip

\textbf{Definition 2.2.} ([CHA; 58]) An abelian monoid $\left( X,\theta
,\oplus \right) $ is called \textit{MV-algebra} if and only if we have an
operation $"^{\prime }"$ such that:

i) $(x^{\prime })^{\prime }=x;$

ii) $x\oplus \theta ^{\prime }=\theta ^{\prime };$

iii) $\left( x^{\prime }\oplus y\right) ^{\prime }\oplus y=$ $\left(
y^{\prime }\oplus x\right) ^{\prime }\oplus x$, for all $x,y\in X.$([Mu;
07]). We denote it by $\left( X,\oplus ,^{\prime },\theta \right) .\medskip $

\textbf{Remark 2.3. }([Mu; 07]) In an MV-algebra the constant element $%
\theta ^{\prime }$ is denoted with $1$, that means

\begin{equation*}
1=\theta ^{\prime },
\end{equation*}%
and the following multiplications are also defined:

\begin{equation*}
x\odot y=\left( x^{\prime }\oplus y^{\prime }\right) ^{\prime },
\end{equation*}%
\begin{equation*}
x\ominus y=x\odot y^{\prime }=\left( x^{\prime }\oplus y\right) ^{\prime }.
\end{equation*}

\textbf{Remark 2.4.} ([COM; 00], Theorem 1.7.1, p. 30)

i) Let $\left( X,\ast ,\theta ,1\right) $ be a bounded commutative
BCK-algebra. If we define 
\begin{equation*}
x^{\prime }=1\ast x,
\end{equation*}%
\begin{equation*}
x\oplus y=1\ast \left( \left( 1\ast x\right) \ast y\right) =\left( x^{\prime
}\ast y\right) ^{\prime },x,y\in X,
\end{equation*}%
we obtain that the algebra $\left( X,\oplus ,^{\prime },\theta \right) $ is
an \textit{MV}-algebra, with 
\begin{equation*}
x\ominus y=x\ast y.
\end{equation*}

ii) If $\left( X,\oplus ,\theta ,^{\prime }\right) $ is an MV-algebra, then $%
\left( X,\ominus ,\theta ,1\right) $ is a bounded commutative
BCK-algebra.\medskip

\textbf{Definition 2.5.}([COM; 00], Definition 4.2.1) An algebra $\left(
W,\circ ,\overline{},1\right) $ of type $\left( 2,1,0\right) ~$is called a 
\textit{Wajsberg algebra (}or\textit{\ W-algebra)} if and only if for every $%
x,y,z\in W$, we have:

i) $1\circ x=x;$

ii) $\left( x\circ y\right) \circ \left[ \left( y\circ z\right) \circ \left(
x\circ z\right) \right] =1;$

iii) $\left( x\circ y\right) \circ y=\left( y\circ x\right) \circ x;$

iv) $\left( \overline{x}\circ \overline{y}\right) \circ \left( y\circ
x\right) =1.\medskip $

\textbf{Remark 2.6. }([COM; 00], Lemma 4.2.2 and Theorem 4.2.5)

i) If $\left( W,\circ ,\overline{},1\right) $ is a Wajsberg algebra,
defining the following multiplications 
\begin{equation*}
x\odot y=\overline{\left( x\circ \overline{y}\right) }
\end{equation*}%
and 
\begin{equation*}
x\oplus y=\overline{x}\circ y,
\end{equation*}%
for all $x,y\in W$, we obtain that $\left( W,\oplus ,\odot ,\overline{}%
,0,1\right) $ is an MV-algebra.

ii) If $\left( X,\oplus ,\odot ,^{\prime },\theta ,1\right) $ is an
MV-algebra, defining on $X$ the operation%
\begin{equation*}
x\circ y=x^{\prime }\oplus y,
\end{equation*}%
it results that $\left( X,\circ ,^{\prime },1\right) $ is a Wajsberg
algebra.\medskip

\textbf{Example 2.7.} We consider the following set $X=\{\theta ,a,b,c,d,e\}$
and the multiplication $"\ast "$ given in the below table:

\begin{center}
\begin{tabular}{l|llllll}
$\ast $ & $\theta $ & $a$ & $b$ & $c$ & $d$ & $e$ \\ \hline
$\theta $ & $\theta $ & $\theta $ & $\theta $ & $\theta $ & $\theta $ & $%
\theta $ \\ 
$a$ & $a$ & $\theta $ & $a$ & $\theta $ & $\theta $ & $\theta $ \\ 
$b$ & $b$ & $b$ & $\theta $ & $b$ & $\theta $ & $\theta $ \\ 
$c$ & $c$ & $a$ & $c$ & $\theta $ & $a$ & $\theta $ \\ 
$d$ & $d$ & $b$ & $a$ & $b$ & $\theta $ & $\theta $ \\ 
$e$ & $e$ & $d$ & $c$ & $b$ & $a$ & $\theta $%
\end{tabular}
\end{center}

Then $\left( X,\ast ,\theta ,1\right) $ becomes a BCK commutative bounded
algebra. The associated MV-algebra is $\left( X,\oplus ,^{\prime },\theta
\right) $, with the multiplication $\oplus $ and the operation $"^{\prime }"$
given in the below tables:

\begin{center}
\begin{tabular}{l|llllll}
$\oplus $ & $\theta $ & $a$ & $b$ & $c$ & $d$ & $e$ \\ \hline
$\theta $ & $\theta $ & $a$ & $b$ & $c$ & $d$ & $e$ \\ 
$a$ & $a$ & $c$ & $d$ & $c$ & $e$ & $e$ \\ 
$b$ & $b$ & $d$ & $b$ & $e$ & $d$ & $e$ \\ 
$c$ & $c$ & $c$ & $e$ & $c$ & $e$ & $e$ \\ 
$d$ & $d$ & $e$ & $d$ & $e$ & $e$ & $e$ \\ 
$e$ & $e$ & $e$ & $e$ & $e$ & $e$ & $e$%
\end{tabular}
~\ 
\begin{tabular}{l|llllll}
$^{\prime }$ & $\theta $ & $a$ & $b$ & $c$ & $d$ & $e$ \\ \hline
& $e$ & $d$ & $c$ & $b$ & $a$ & $\theta $%
\end{tabular}%
.
\end{center}

([WDH; 17], Example 3.3).

The associated Wajsberg algebra $\left( X,\circ ,^{\prime },1\right) $ is:

\begin{center}
\begin{tabular}{l|llllll}
$\circ $ & $\theta $ & $a$ & $b$ & $c$ & $d$ & $e$ \\ \hline
$\theta $ & $e$ & $e$ & $e$ & $e$ & $e$ & $e$ \\ 
$a$ & $d$ & $e$ & $d$ & $e$ & $e$ & $e$ \\ 
$b$ & $c$ & $c$ & $e$ & $c$ & $e$ & $e$ \\ 
$c$ & $b$ & $d$ & $b$ & $e$ & $d$ & $e$ \\ 
$d$ & $a$ & $c$ & $d$ & $c$ & $e$ & $e$ \\ 
$e$ & $\theta $ & $a$ & $b$ & $c$ & $d$ & $e$%
\end{tabular}
\end{center}

\textbf{Proposition 2.8. }([Bu; 06]) \textit{Let} $\left( X,\oplus ,^{\prime
},\theta \right) $ \textit{be an MV-algebra. We have that} 
\begin{equation*}
x\oplus x^{\prime }=1.
\end{equation*}

\textbf{Proof.} Indeed, from Definition 2.2 ii) and iii), it results that 
\newline
$1=\left( x^{\prime }\oplus 1\right) ^{\prime }\oplus 1=\left( 1^{\prime
}\oplus x\right) ^{\prime }\oplus x=x^{\prime }\oplus x._{{}}\medskip $

\textbf{Proposition 2.9. }([Mu; 07]) \textit{Let} $\left( X,\oplus ,^{\prime
},\theta \right) $ \textit{be an MV-algebra. For }$x,y\in X$\textit{, the
following statements are equivalent:}

\textit{i)} $x^{\prime }\oplus y=1;$

\textit{ii)} $x\odot y^{\prime }=0;$

\textit{iii)} $y=x\oplus \left( y\ominus x\right) =x\oplus \left( y^{\prime
}\oplus x\right) ^{\prime };$

\textit{iv)} \textit{There is an element} $z\in X$ \textit{such that}~$%
x\oplus z=y$.\medskip $_{{}}\medskip $

\textbf{Definition 2.10.} ([COM; 00]) Let $\left( X,\oplus ,^{\prime
},\theta \right) $ be an MV-algebra and $x,y\in X$. On $X$, we define the
following order relation:%
\begin{equation*}
x\leq y\text{ if and only if }x^{\prime }\oplus y=1.
\end{equation*}%
$\medskip $

\textbf{Remark} \textbf{2.11.} It is clear that $x\leq y$ if and only if $x$
and $y$ satisfy one of the equivalent conditions i)-iv) from Proposition 2.9.%
\textbf{\medskip }

\textbf{Definition 2.12.} [FRT; 84] If $\left( W,\circ ,\overline{},1\right) 
$ is a Wajsberg algebra, on $W$ we define the following binary relation 
\begin{equation}
x\leq y~\text{if~and~only~if~}x\circ y=1.  \tag{2.1.}
\end{equation}%
This relation is an order relation, called \textit{the natural order
relation on }$W$\textit{.}

For other details regarding MV-algebras and Wajsberg algebras, the readers
are referred to [Io; 08] and [Pi; 07].

\begin{equation*}
\end{equation*}

\textbf{3. Block codes associated to MV-algebras and Wajsberg algebras }%
\begin{equation*}
\end{equation*}

In [JUN; 11], were introduced binary block-codes over finite BCK-algebras.
In a similar way, code over MV-algebras and Wajsberg algebras can be
established.

Let $A$ be a nonempty set and let $\left( X,\oplus ,^{\prime },\theta
\right) $ be an MV-algebra.\medskip

\textbf{Definition 3.1.} A mapping $f:A\rightarrow X$ is called an \textit{%
MV-function} on $A.$ A \textit{cut function} \textit{of} $f$ is a map $%
f_{r}:A\rightarrow \{0,1\},r\in X$, such that 
\begin{equation*}
f_{r}\left( x\right) =1\text{, if and only if \ }r^{\prime }\oplus f\left(
x\right) =1\text{,}\forall x\in A\text{.}
\end{equation*}%
\newline
A \textit{cut subset} of $A$ is the following subset of $A$

\begin{equation*}
A_{r}=\{x\in A~/~r^{\prime }\oplus f\left( x\right) =1\}.
\end{equation*}

\textbf{Remark 3.2.} If $y\in A_{r}\cup A_{s},$therefore $y\in A_{r\odot s}.$
Indeed, for $y\in A_{r}\cup A_{s}$, we suppose that $y\in A_{s}$. From here,
we have $y\in A_{r\odot s}$, since $A_{r\odot s}=\{x\in A~/~\left( r\odot
s\right) ^{\prime }\oplus f\left( x\right) =1\}$ and $\left( r\odot s\right)
^{\prime }\oplus f\left( y\right) =(r^{\prime }\oplus s^{\prime })\oplus
f\left( y\right) =r^{\prime }\oplus (s^{\prime }\oplus f\left( y\right)
)=r^{\prime }\oplus 1=1.\medskip $

\textbf{Remark 3.3.} Let $f:A\rightarrow X$ \ be an MV-function on $A.$We
define on $X$ the following binary relation 
\begin{equation*}
\forall r,s\in X,r\sim s~~\text{if~and~only~if~~}A_{r}=A_{s}.
\end{equation*}%
This relation is an equivalence relation on $X$ and we denote with $%
\widetilde{r}$ the equivalence class of an element $r\in X.\medskip $

\bigskip \textbf{Proposition 3.4}. \textit{Let} $f:A\rightarrow X$ \ \textit{%
be an MV-function on} $A$. \textit{Therefore}%
\begin{equation*}
f\left( x\right) =\sup \{r\in X~/~f_{r}\left( x\right) =1\}.
\end{equation*}

\textbf{Proof.} For a chosen element $x\in A$, we denote $f\left( x\right)
=s,s\in X$. From here we have $s^{\prime }\oplus f\left( x\right) =s^{\prime
}\oplus s=1$, therefore $f_{s}\left( x\right) =1.$ If there is an element $%
r\in X$ such that $r^{\prime }\oplus f\left( x\right) =r^{\prime }\oplus s=1$%
, we obtain that $r\leq s$. We consider the set $M=\{r\in X~/~f_{r}\left(
x\right) =1\}$. Since $f_{s}\left( x\right) =1$, it results that $s\in M$,
therefore $f\left( x\right) =s=\sup M.\Box \medskip $

The above proposition extends to MV-algebras results obtained in Proposition
3.4 from [JUN; 11]\medskip

\textbf{Proposition 3.5.} \textit{Let} $f:A\rightarrow X$ \ \textit{be an
MV-function on} $A$. \textit{Therefore, for} $r,s\in X$, \textit{we have
that } $r\leq s$ \textit{implies} $A_{s}\subseteq A_{r}$. \textit{If} $A=X$ 
\textit{and} $f:X\rightarrow X$ \textit{is the identity function,} $f\left(
x\right) =x$, \textit{then the converse of this statement is also
true.\medskip }

\textbf{Proof.} Let $r,s\in X$ such that $r\leq s$. From Proposition 2.9.,
iv), it results that $s=r\oplus t,t\in X$. For $x\in A_{s}$, we have $%
s^{\prime }\oplus f\left( x\right) =1$, which is equivalent with $s\leq
f\left( x\right) $, that means $f\left( x\right) =s\oplus q,q\in X$. It
results that $f\left( x\right) =\left( r\oplus t\right) \oplus q=r\oplus
(t\oplus q)$, therefore $r\leq f\left( x\right) .$ From here, we obtain that 
$r^{\prime }\oplus f\left( x\right) =1$, then $x\in A_{r}$.

For the converse, if $A_{s}\subseteq A_{r}$, we have that $s\in
A_{s}\subseteq A_{r}$, therefore $s\in A_{r}$. It results that $r^{\prime
}\oplus f\left( s\right) =r^{\prime }\oplus s=1$, which implies $r\leq s$.$%
\Box \medskip $

The above proposition extends to MV-algebras results obtained in Proposition
3.6 from [JUN; 11]\medskip .\medskip

Let $A$ be a set with $n$ elements. We consider $A=\{1,2,...,n\}$ and let $X$
be an MV-algebra. Using above notations, to each equivalence class $%
\widetilde{x},x\in X,$ will correspond the codeword $%
w_{x}=f_{x}=x_{1}x_{2}...x_{n}$, with $x_{i}=j,$ if and only if $f_{x}\left(
i\right) =j,i\in A,j\in \{0,1\}$.We denote this code with $V_{X}$. In this
way, each MV-function $f:A\rightarrow X$ has associated a binary block-code
of length $n$ and $V_{X}=\{f_{x},x\in X\}$ $.$

Let $V$ be a binary block-code of length $n$ and $w_{x}=x_{1}x_{2}...x_{n}%
\in V,$ $w_{y}=y_{1}y_{2}...y_{n}\in V$ ~be two codewords$.$ On $V$ we can
define the following partial order relation: 
\begin{equation}
w_{x}\preceq w_{y}\text{ if and only if }y_{i}\leq x_{i},i\in \{1,2,...,n\}.
\tag{3.1.}
\end{equation}

\textbf{Proposition 3.6.} \textit{Let} $X$ \textit{be an MV-algebra. With
the above notations, relation} $A_{s}\subseteq A_{r}$ \textit{is equivalent
with} $f_{r}\preceq f_{s}$.\medskip

\textbf{Proof.} Assuming $A_{s}\subseteq A_{r}$, we have $f_{s}\left(
x\right) =f_{r}\left( x\right) =1$, for all $x\in A_{s}$. It results $%
f_{s}\left( x\right) \leq f_{r}\left( x\right) $, for all $x\in X$, that
means $f_{r}\preceq f_{s}$.

Conversely, if $f_{r}\preceq f_{s}$, we have $f_{s}\left( x\right) \leq
f_{r}\left( x\right) $, for all $x\in X$. If $x\in A_{s}$, we obtain $%
s^{\prime }\oplus f\left( x\right) =1$, which is equivalent with $%
f_{s}\left( x\right) =1$. This implies that $f_{r}\left( x\right) =1$, which
is equivalent with $r^{\prime }\oplus f\left( x\right) =1$. From here, we
obtain $x\in A_{r}$.$\Box \medskip $

\textbf{Proposition} \textbf{3.7}. \textit{Let} $\left( X,\oplus ,^{\prime
},\theta \right) $ \ \textit{be a finite MV-algebra.} \textit{To algebra} $X$
\textit{corresponds a block-code} $V$ \textit{such that} $\left( X,\leq
\right) $ \textit{is isomorphic to} $\left( V_{X},\preceq \right) $ \textit{%
as ordered sets.}\medskip 

\textbf{Proof.} Let $\left( X,\oplus ,^{\prime },\theta \right) $ be a
finite MV-algebra and, for $A=X$, let $f:X\rightarrow X$ be the identity
function which is an MV-function. The function $f$ generates the following
set $M=\{f_{r}\ /\ r\in X\}$ of cuts functions.~Then the set $M$ with the
order $\preceq $ is the associated block-code to the MV-algebra $X$, with
the order relation defined in $\left( 3.1\right) $. Let $g:X\rightarrow M,$ $%
g\left( r\right) =f_{r},$ for all $r\in X$. We will prove that this map is
bijective and $r\leq s$ is equivalent with $f_{r}\preceq f_{s}$. By
definition, the map is surjective. The map $g$ is injective. Indeed, if $%
r,s\in X$ with $g\left( s\right) =g\left( r\right) ,$ we have $f_{r}=f_{s}$,
therefore $A_{r}=A_{s}$. It results that $1=r^{\prime }\oplus f\left(
s\right) =r^{\prime }\oplus s$ and $1=s^{\prime }\oplus f\left( r\right)
=s^{\prime }\oplus r$, therefore, $r\leq s$ and $s\leq r$. From here, we
obtain that $r=s$, therefore $g$ is a bijective map. Let $r,s\in X$ such
that $r\leq s$. From Proposition 3.5, this is equivalent with $%
A_{s}\subseteq A_{r}$, which is equivalent with $f_{r}\preceq f_{s}$, from
Proposition 3.6. $\Box \medskip $

\textbf{Theorem 3.8}. \textit{Let} $A=\left( X,\ast ,\theta ,1\right) $ 
\textit{be a bounded commutative BCK-algebra and} $B=\left( X,\oplus
,^{\prime },\theta \right) $ \textit{be} \textit{the associated MV-algebra.} 
\textit{Therefore} $A$ \textit{and} $B$ \textit{are code equivalent, that
means determine the same binary block-code.}\medskip

\textbf{Proof.} With the above notations, let $V_{A}$ be the code associated
to the BCK algebra $A$ and $V_{B}$ be the code associated to the MV-algebra $%
B$. Let $\varphi :V_{A}\rightarrow V_{B},\varphi \left( f_{r}\right)
=f_{r^{\prime }}^{\prime }$, where $f_{r}$ is a cut function associated to
BCK-algebra as in [JUN; 11], $f_{r^{\prime }}^{\prime }~$is the corresponded
cut function associated to MV-algebra and $r^{\prime }$ is the complement of
the element $r\in X$. The map $\varphi $ is injective, therefore bijective.
Indeed, if $\varphi \left( f_{r}\right) =\varphi \left( f_{s}\right) $, it
results $f_{r^{\prime }}^{\prime }=f_{s^{\prime }}^{\prime },$ therefore $%
r^{\prime }=s^{\prime }$ and $r=s.\Box \medskip $

\textbf{Definition 3.9.} ([CHA; 19]) Let $\left( X,\oplus ,^{\prime },\theta
\right) $ \ be an MV-algebra. The distance function defined on the algebra $%
A $ is:%
\begin{eqnarray*}
d &:&X\times X\rightarrow X,d\left( x,y\right) =\left( x\ominus y\right)
\oplus \left( y\ominus x\right) = \\
&=&(x\odot y^{\prime })\oplus \left( y\odot x^{\prime }\right) = \\
&=&\left( x^{\prime }\oplus y\right) ^{\prime }\oplus \left( y^{\prime
}\oplus x\right) ^{\prime }.
\end{eqnarray*}

In the following, inspired from the above definition, we will define a new
distance on a finite MV-algebra $X$ with $n$ elements. Let $\left( X,\oplus
,^{\prime },\theta \right) $ \ be an MV-algebra, $f:A\rightarrow X$ be an 
\textit{MV-function} on $A$, $f_{r}:A\rightarrow \{0,1\},~r\in X,$ be a cut
function and $A_{r}$ be the associated cut subset. Let $\varphi
:X\rightarrow \mathcal{P}\left( A\right) ,\varphi \left( r\right) =A_{r}$.
Since $\left( \mathcal{P}\left( A\right) ,\cup ,\cap ,^{\prime },\emptyset
,A\right) $ is an MV-algebra, let $d$ the distance on $\mathcal{P}\left(
A\right) $ defined as above. We define on $X$ the following map 
\begin{equation*}
D:X\times X\rightarrow \mathbb{R}_{+},
\end{equation*}%
\begin{equation*}
D\left( r,s\right) =\left\vert d\left( \varphi \left( r\right) ,\varphi
\left( s\right) \right) \right\vert
\end{equation*}

\textbf{Proposition 3.10.}

\textit{1)} $D\left( r,s\right) =0$ \textit{if and only if} $r=s$.

\textit{2)} $D\left( r,s\right) \leq D\left( r,t\right) +D\left( t,s\right)
. $

\textit{3)} $D\left( r,\theta \right) =\left\vert A_{r}\right\vert $ \textit{%
and} $D\left( r,1\right) =\left\vert A_{r}^{\prime }\right\vert $, where $%
A_{r}^{\prime }$ is the complement of the set $A_{r}^{\prime }.\medskip $

\textbf{Proof.}

1) Indeed, if $D\left( r,s\right) =\left\vert d\left( \varphi \left(
r\right) ,\varphi \left( s\right) \right) \right\vert =0$, we have $%
(A_{r}^{\prime }\cap A_{s})\cup \left( A_{r}\cap A_{s}^{\prime }\right)
=\emptyset .$ From here, we have that $A_{r}^{\prime }\cap A_{s}=A_{r}\cap
A_{s}^{\prime }=\emptyset $, therefore $A_{r}=A_{s}$ which is equivalent
with $r=s.$

2) Indeed, $D\left( r,s\right) =\left\vert d\left( \varphi \left( r\right)
,\varphi \left( s\right) \right) \right\vert \leq $\newline
$\leq \left\vert d\left( \varphi \left( r\right) ,\varphi \left( t\right)
\right) +d\left( \varphi \left( t\right) ,\varphi \left( s\right) \right)
\right\vert \leq \left\vert d\left( \varphi \left( r\right) ,\varphi \left(
t\right) \right) \right\vert +\left\vert d\left( \varphi \left( t\right)
,\varphi \left( s\right) \right) \right\vert .\Box \medskip $

\textbf{Definition 3.11.} The distance $D$ is called \textit{the Hamming
distance between the elements} $r$ and $r,s\in X,$ and it is the really
Hamming distance between their associated code-words $f_{r}$ and $f_{s}$%
.\medskip 

In a similar way as above, we can introduce binary block-codes over finite
Wajsberg algebra.

Let $A$ be a nonempty set and let $\left( W,\circ ,\overline{},1\right) $ be
a Wajsberg algebra.\medskip

\textbf{Definition 3.12.} A mapping $f:A\rightarrow W$ is called a \textit{%
W-function} on $A$. A \textit{cut function} \textit{of} $f$ is a map $%
f_{r}:A\rightarrow \{0,1\},r\in W$, such that 
\begin{equation*}
f_{r}\left( x\right) =1\text{, if and only if \ }r\circ f\left( x\right) =1%
\text{,}\forall x\in A\text{.}
\end{equation*}%
\newline
A \textit{cut subset} of $A$ is the following subset of $A$

\begin{equation*}
A_{r}=\{x\in A~/~r\circ f\left( x\right) =1\}.
\end{equation*}%
Let $f:A\rightarrow W$ \ be an W-function on $A.$We define on $W$ the
following binary relation 
\begin{equation*}
\forall r,s\in W,r\sim s~~\text{if~and~only~if~~}A_{r}=A_{s}.
\end{equation*}%
This relation is an equivalence relation on $W$ and we denote with $%
\widetilde{r}$ the equivalence class of an element $r\in W$.\medskip 

\textbf{Proposition 3.13.} \textit{Let} $f:A\rightarrow W$ \ \textit{be a
W-function on} $A$. \textit{Therefore, for} $r,s\in W$, \textit{we have that 
} $r\leq s$ \textit{implies} $A_{s}\subseteq A_{r}$. \textit{If} $A=W$ 
\textit{and} $f:W\rightarrow W$ \textit{is the identity function,} $f\left(
x\right) =x$, \textit{then the converse of this statement is also true}%
.\medskip

\textbf{Proof.} Assuming that $r\leq s$, we have $r\circ s=1$. Let $w\in
A_{s}$, therefore $s\circ f\left( w\right) =1$. We will prove that $r\circ
f\left( w\right) =1$. Since for all $x,y,z\in W$, from $\left( x\circ
y\right) \circ \left[ \left( y\circ z\right) \circ \left( x\circ z\right) %
\right] =1$, we have $\left( r\circ s\right) \circ \left[ \left( s\circ
f\left( w\right) \right) \circ \left( r\circ f\left( w\right) \right) \right]
=1$, therefore $1\circ \left[ 1\circ \left( r\circ f\left( w\right) \right) %
\right] =1$. From here, since $1\circ x=x$, it results $1=r\circ f\left(
w\right) $, therefore $A_{s}\subseteq A_{r}$.

For the converse, if $A_{s}\subseteq A_{r}$, we have that $s\in A_{r}$,
therefore $1=r\circ f\left( s\right) =r\circ s$. We obtain $r\leq s$.$\Box
\medskip $

The above proposition extends to W-algebras results obtained in Proposition
3.6 from [JUN; 11]\medskip .\medskip

Let $A$ be a set with $n$ elements. We consider $A=\{1,2,...,n\}$ and let $W$
be a W-algebra. Using above notations, to each equivalence class $\widetilde{%
x},x\in W$, will correspond the codeword $w_{x}=f_{x}=x_{1}x_{2}...x_{n}$,
with $x_{i}=j$, if and only if $f_{x}\left( i\right) =j,i\in A,j\in \{0,1\}$%
.We denote this code with $V_{W}$. In this way, each W-function $%
f:A\rightarrow W$ \ has associated a binary block-code of length $n$ and $%
V_{W}=\{f_{x},x\in W\}$.\medskip

\textbf{Proposition 3.14.} \textit{Let} $W$ \textit{be a W-algebra. With the
above notations, relation} $A_{s}\subseteq A_{r}$ \textit{is equivalent with}
$f_{r}\preceq f_{s}$.\medskip

\textbf{Proof.} By straightforward calculation. $\Box \medskip $

\textbf{Proposition} \textbf{3.15}. \textit{Let} $\left( W,\circ ,\overline{}%
,1\right) $ \ \textit{be a finite W-algebra.} \textit{To algebra} $W$ 
\textit{corresponds a block-code} $V$ \textit{such that} $\left( W,\leq
\right) $ \textit{is isomorphic to} $\left( V_{W},\preceq \right) $ \textit{%
as ordered sets.}\medskip 

\textbf{Proof.} By straightforward calculation. $\Box \medskip $

\textbf{Theorem 3.16}. \textit{Let} $A=\left( X,\oplus ,^{\prime },\theta
\right) $ \textit{be an MV-algebra and} $B=\left( X,\circ ,\overline{}%
,1\right) $ \textit{be the associated Wajsberg algebra. Therefore} $A$ 
\textit{and} $B$ \textit{are code equivalent, that means determine the same
binary block-code.}\medskip

\textbf{Proof.} With the above notations, let $V_{A}$ be the code associated
to the MV-algebra $A$ and $V_{B}$ be the code associated to the W-algebra $B$%
. Let $\varphi :V_{A}\rightarrow V_{B},\varphi \left( f_{r}\right)
=f_{r}^{\prime }$, where $f_{r}$ is a cut function associated to MV-algebra
and $f_{r}^{\prime }$ is the corresponded cut function associated to
W-algebra. The map $\varphi $ is injective, therefore bijective. Indeed, if $%
\varphi \left( f_{r}\right) =\varphi \left( f_{s}\right) $, it results $%
f_{r}^{\prime }=f_{s}^{\prime },$ therefore $r=s$. $\Box \medskip $

\textbf{Theorem 3.17}. \textit{Let} $A=\left( X,\ast ,\theta ,1\right) $ 
\textit{be a bounded commutative BCK-algebra,} $B=\left( X,\oplus ,^{\prime
},\theta \right) $ \textit{be} \textit{the associated MV-algebra and }$%
C=\left( X,\circ ,\overline{},1\right) $ \textit{be the associated Wajsberg
algebra.} \textit{Therefore} $A,B$ \textit{and} $C$ \textit{are code
equivalent, that means determine the same binary block-code.}\medskip
\medskip

\textbf{Proof.} The result is obtained from Theorem 3.8 and Theorem 3.16. $%
\Box \medskip $

\textbf{Example 3.18.} Using algebras from Example 2.7\textbf{,} we can see
that the code associated to BCK commutative bounded algebra $\left( X,\ast
,\theta ,1\right) $ is\newline
$V_{1}=\{111111,010111,001011,000101,000011,000001\}$, the code associated
to MV-algebra $\left( X,\oplus ,^{\prime },\theta \right) $,$~V_{2}$, is the
same with $V_{1}$ and the code associated to Wajsberg algebra $\left(
X,\circ ,\overline{},1\right) $,$~V_{3}$, is also the same with $V_{1}$.

i) The BCK commutative bounded algebra and the associated code are given in
the below tables:

\ \ {\small \qquad \qquad\ \ \ \ \ }

\begin{equation*}
\begin{tabular}{l|llllll}
$\ast $ & $\theta $ & $a$ & $b$ & $c$ & $d$ & $e$ \\ \hline
$\theta $ & $\theta $ & $\theta $ & $\theta $ & $\theta $ & $\theta $ & $%
\theta $ \\ 
$a$ & $a$ & $\theta $ & $a$ & $\theta $ & $\theta $ & $\theta $ \\ 
$b$ & $b$ & $b$ & $\theta $ & $b$ & $\theta $ & $\theta $ \\ 
$c$ & $c$ & $a$ & $c$ & $\theta $ & $a$ & $\theta $ \\ 
$d$ & $d$ & $b$ & $a$ & $b$ & $\theta $ & $\theta $ \\ 
$e$ & $e$ & $d$ & $c$ & $b$ & $a$ & $\theta $%
\end{tabular}%
~~~%
\begin{tabular}{l|llllll}
$\ast $ & $\theta $ & $a$ & $b$ & $c$ & $d$ & $e$ \\ \hline
$\theta $ & $\mathbf{1}$ & $\mathbf{1}$ & $\mathbf{1}$ & $\mathbf{1}$ & $%
\mathbf{1}$ & $\mathbf{1}$ \\ 
$a$ & $0$ & $\mathbf{1}$ & $0$ & $\mathbf{1}$ & $\mathbf{1}$ & $\mathbf{1}$
\\ 
$b$ & $0$ & $0$ & $\mathbf{1}$ & $0$ & $\mathbf{1}$ & $\mathbf{1}$ \\ 
$c$ & $0$ & $0$ & $0$ & $\mathbf{1}$ & $0$ & $\mathbf{1}$ \\ 
$d$ & $0$ & $0$ & $0$ & $0$ & $\mathbf{1}$ & $\mathbf{1}$ \\ 
$e$ & $0$ & $0$ & $0$ & $0$ & $0$ & $\mathbf{1}$%
\end{tabular}%
.
\end{equation*}%
\medskip

ii) The MV-algebra algebra and the associated code are given in the below
tables:

\begin{equation*}
\begin{tabular}{l|llllll}
$\oplus $ & $\theta $ & $a$ & $b$ & $c$ & $d$ & $e$ \\ \hline
$\theta $ & $\theta $ & $a$ & $b$ & $c$ & $d$ & $e$ \\ 
$a$ & $a$ & $c$ & $d$ & $c$ & $e$ & $e$ \\ 
$b$ & $b$ & $d$ & $b$ & $e$ & $d$ & $e$ \\ 
$c$ & $c$ & $c$ & $e$ & $c$ & $e$ & $e$ \\ 
$d$ & $d$ & $e$ & $d$ & $e$ & $e$ & $e$ \\ 
$e$ & $e$ & $e$ & $e$ & $e$ & $e$ & $e$%
\end{tabular}%
~~~%
\begin{tabular}{l|llllll}
$\oplus $ & $\theta $ & $a$ & $b$ & $c$ & $d$ & $e$ \\ \hline
$\theta $ & $0$ & $0$ & $0$ & $0$ & $0$ & $\mathbf{1}$ \\ 
$a$ & $0$ & $0$ & $0$ & $0$ & $\mathbf{1}$ & $\mathbf{1}$ \\ 
$b$ & $0$ & $0$ & $0$ & $\mathbf{1}$ & $0$ & $\mathbf{1}$ \\ 
$c$ & $0$ & $0$ & $\mathbf{1}$ & $0$ & $\mathbf{1}$ & $\mathbf{1}$ \\ 
$d$ & $0$ & $\mathbf{1}$ & $0$ & $\mathbf{1}$ & $\mathbf{1}$ & $\mathbf{1}$
\\ 
$e$ & $\mathbf{1}$ & $\mathbf{1}$ & $\mathbf{1}$ & $\mathbf{1}$ & $\mathbf{1}
$ & $\mathbf{1}$%
\end{tabular}%
\end{equation*}%
\medskip

iii) \ The Wajsberg algebra and the associated code are given in the below
tables:

~%
\begin{equation*}
\begin{tabular}{l|llllll}
$\circ $ & $\theta $ & $a$ & $b$ & $c$ & $d$ & $e$ \\ \hline
$\theta $ & $e$ & $e$ & $e$ & $e$ & $e$ & $e$ \\ 
$a$ & $d$ & $e$ & $d$ & $e$ & $e$ & $e$ \\ 
$b$ & $c$ & $c$ & $e$ & $c$ & $e$ & $e$ \\ 
$c$ & $b$ & $d$ & $b$ & $e$ & $d$ & $e$ \\ 
$d$ & $a$ & $c$ & $d$ & $c$ & $e$ & $e$ \\ 
$e$ & $\theta $ & $a$ & $b$ & $c$ & $d$ & $e$%
\end{tabular}%
\medskip \ ~~%
\begin{tabular}{l|llllll}
$\circ $ & $\theta $ & $a$ & $b$ & $c$ & $d$ & $e$ \\ \hline
$\theta $ & $\mathbf{1}$ & $\mathbf{1}$ & $\mathbf{1}$ & $\mathbf{1}$ & $%
\mathbf{1}$ & $\mathbf{1}$ \\ 
$a$ & $0$ & $\mathbf{1}$ & $0$ & $\mathbf{1}$ & $\mathbf{1}$ & $\mathbf{1}$
\\ 
$b$ & $0$ & $0$ & $\mathbf{1}$ & $0$ & $\mathbf{1}$ & $\mathbf{1}$ \\ 
$c$ & $0$ & $0$ & $0$ & $\mathbf{1}$ & $0$ & $\mathbf{1}$ \\ 
$d$ & $0$ & $0$ & $0$ & $0$ & $\mathbf{1}$ & $\mathbf{1}$ \\ 
$e$ & $0$ & $0$ & $0$ & $0$ & $0$ & $\mathbf{1}$%
\end{tabular}%
\end{equation*}

We remark that the attached code for each of these algebras is as a skeleton
for these algebras( the same for all three algebras), on which we can insert
different structures, as for example: a BCK commutative bounded algebra or
an MV-algebra or a Wajsberg algebra. We consider the algebra \thinspace $%
A\in $\{BCK, MV, Wajsberg\}, finite with $n$ elements. The \ skeleton of
such an algebra is\ a matrix of order $n$ in which the elements of this
matrix is black or white squares, with black square on the position $\left(
i,j\right) $, if and only if $x_{i}\leq x_{j}$ in $A$, $x_{i},x_{j}\in A$.
The associated skeleton of such a finite algebra \thinspace $A\in $\{BCK,
MV, Wajsberg\} is nothing else than a representation of the associated order
relation on the algebra $A$.

As we can see in the below two tables, the skeleton is the same for the BCK
bounded commutative algebra and for the attached Wajsberg algebra. For the
attached MV algebra, the skeleton is the same but symmetric in respect to
the Ox axis. The skeleton generates the same order relation on $A$ as the
attached binary block code on $A,V_{A}$.\medskip

\begin{equation}
\begin{tabular}{l|llllll}
$\ast ,\circ $ & $\theta $ & $a$ & $b$ & $c$ & $d$ & $e$ \\ \hline
$\theta $ & $\mathbf{\blacksquare }$ & $\mathbf{\blacksquare }$ & $\mathbf{%
\blacksquare }$ & $\mathbf{\blacksquare }$ & $\mathbf{\blacksquare }$ & $%
\mathbf{\blacksquare }$ \\ 
$a$ &  & $\mathbf{\blacksquare }$ &  & $\mathbf{\blacksquare }$ & $\mathbf{%
\blacksquare }$ & $\mathbf{\blacksquare }$ \\ 
$b$ &  &  & $\mathbf{\blacksquare }$ &  & $\mathbf{\blacksquare }$ & $%
\mathbf{\blacksquare }$ \\ 
$c$ &  &  &  & $\mathbf{\blacksquare }$ &  & $\mathbf{\blacksquare }$ \\ 
$d$ &  &  &  &  & $\mathbf{\blacksquare }$ & $\mathbf{\blacksquare }$ \\ 
$e$ &  &  &  &  &  & $\mathbf{\blacksquare }$%
\end{tabular}%
~%
\begin{tabular}{l|llllll}
$\oplus $ & $\theta $ & $a$ & $b$ & $c$ & $d$ & $e$ \\ \hline
$\theta $ &  &  &  &  &  & $\mathbf{\blacksquare }$ \\ 
$a$ &  &  &  &  & $\mathbf{\blacksquare }$ & $\mathbf{\blacksquare }$ \\ 
$b$ &  &  &  & $\mathbf{\blacksquare }$ &  & $\mathbf{\blacksquare }$ \\ 
$c$ &  &  & $\mathbf{\blacksquare }$ &  & $\mathbf{\blacksquare }$ & $%
\mathbf{\blacksquare }$ \\ 
$d$ &  & $\mathbf{\blacksquare }$ &  & $\mathbf{\blacksquare }$ & $\mathbf{%
\blacksquare }$ & $\mathbf{\blacksquare }$ \\ 
$e$ & $\mathbf{\blacksquare }$ & $\mathbf{\blacksquare }$ & $\mathbf{%
\blacksquare }$ & $\mathbf{\blacksquare }$ & $\mathbf{\blacksquare }$ & $%
\mathbf{\blacksquare }$%
\end{tabular}%
\medskip  \tag{3.2.}
\end{equation}

\begin{equation*}
\end{equation*}

\textbf{4. Some remarks regarding Wajsberg algebras}%
\begin{equation*}
\end{equation*}

\textbf{Remark 4.1.} Let\textbf{\ }$\left( X,\leq \right) $\textbf{\ }be a
finite totally ordered set, $X=\{x_{0},x_{1},...,x_{n}\}$, with $x_{0}$ the
first element and $x_{n}$ the last element. Using this order relation, we
define the following multiplication $"\circ "$ on $X$:%
\begin{equation}
\left\{ 
\begin{array}{c}
x_{i}\circ x_{j}=1\text{, if }x_{i}\leq x_{j}; \\ 
x_{i}\circ x_{j}=x_{n-i+j}\text{, otherwise;} \\ 
x_{0}=\theta ,x_{n}=1,x\circ \theta =\overline{x}.%
\end{array}%
\right.  \tag{4.1.}
\end{equation}%
Therefore, $\left( X,\circ ,\overline{},1\right) $ is a Wajsberg algebra. We
remark that this is the only way to define a W-algebra structure on a finite
totally ordered set such that the induced order relation on this algebra is
given in $(2.1)$. We also remark that $\overline{x}_{i}=x_{n-1}$. ([FRT;
84], Theorem 19).\medskip

\textbf{Definition 4.2.} Let $\left( W_{1},\circ ,\overline{},1\right) $ and 
$\left( W_{2},\cdot ,^{\prime },1\right) $ be two finite Wajsberg algebras.
We define on the Cartesian product of these algebras, $W=W_{1}\times W_{2}$,
the following multiplication $"\nabla "$,%
\begin{equation}
\left( x_{1},x_{2}\right) \nabla \left( y_{1},y_{2}\right) =\left(
x_{1}\circ y_{1},x_{2}\cdot y_{2}\right) ,  \tag{4.2.}
\end{equation}%
The complement of the element $\left( x_{1},x_{2}\right) $ is $\rceil \left(
x_{1},x_{2}\right) =$ $\left( \overline{x}_{1},x_{2}^{\prime }\right) $ and $%
\mathbf{1}=\left( 1,1\right) $. Therefore, by straightforward calculation,
we obtain that $\left( W,\nabla ,\rceil ,\mathbf{1}\right) $ is also a
Wajsberg algebra.\medskip

\textbf{Remark 4.3.} If $x=\left( x_{1},x_{2}\right) ,y=\left(
y_{1},y_{2}\right) \in W$, then the order relation corresponded to the
algebra $\left( W,\nabla ,\rceil ,\mathbf{1}\right) $ is given as follow:%
\begin{equation}
x\leq _{W}y\text{ if and only if }x_{1}\leq _{W_{1}}y_{1}\text{ and }%
x_{2}\leq _{W_{2}}y_{2}\text{.}  \tag{4.3.}
\end{equation}

\textbf{Definition 4.4.} Let $\left( W_{1},\circ ,\overline{},1\right) $ and 
$\left( W_{2},\cdot ,^{\prime },1\right) $ be two Wajsberg algebras. A map $%
f:W_{1}\rightarrow W_{2}$ is a morphism of Wajsberg algebras if and only if:

1) $f\left( 0\right) =0;$

2) $f\left( x\circ y\right) =f\left( x\right) \cdot f\left( y\right) ;$

3) $f\left( \overline{x}\right) =\left( f\left( x\right) \right) ^{\prime
}.\medskip $

\textbf{Proposition 4.5.} \textit{The algebras} $W_{1}\times W_{2}$ \textit{%
and} $W_{2}\times W_{1\text{ }}$\textit{are isomorphic.\medskip }

\textbf{Proof.} Let $f:$ $W_{1}\times W_{2}\rightarrow W_{2}\times W_{1}$, $%
f\left( \left( x_{1},x_{2}\right) \right) =\left( x_{2},x_{1}\right) $ is a
bijective morphism. $\Box \medskip $

\textbf{Remark 4.6.} 1) ([HR; 99], Theorem 5.2, p. 43) An MV-algebra is
finite if and only if it is isomorphic to a finite product of totally
ordered MV algebras.

2) If \ $M=\left( X,\oplus ,\odot ,^{\prime },\theta ,1\right) $ is a
totally ordered MV-algebra, then the obtained Wasjberg algebra, $\,W=\left(
X,\circ ,^{\prime },1\right) $, is also totally ordered. The converse is
also true. Indeed, since $x\circ y=x^{\prime }\oplus y$, we have that $x\leq
_{M}y$ if and only $x\leq _{W}y$.

3) Using connections between MV-algebras and Wajsberg algebras, from the
above, we have that a Wajsberg algebra is finite if and only if it is
isomorphic to a finite product of totally ordered Wajsberg algebras.

4) If an MV-algebra or a Wajsberg algebra are finite with a prime number of
elements, therefore these algebras are totally ordered.

5) If two Wajsberg algebras are isomorphic, then these algebras are also
isomorphic as ordered sets.\medskip

\textbf{Definition 4.7.} Let $n$ be a natural number, $n\geq 4$. We consider
the decomposition of the number $n$ in factors: 
\begin{equation*}
n=q_{1}q_{2}...q_{t},q_{i}\in N,1<q_{i}<n,i\in \{1,2,...,t\}.
\end{equation*}%
This decomposition is not unique. We will count only one time the
decompositions with the same terms but with other order of them in the
product. We denote with $\pi _{n}$ the number of all such
decompositions.\medskip

From the above, we obtain the following Theorem.\medskip

\textbf{Theorem 4.8.} \textit{Let} $n$ \textit{be a natural number,} $n\geq 2
$. \textit{There are exactly} $\pi _{n}$ \textit{nonismorphic(as ordered
sets) Wajsberg algebras with} $n$ \textit{elements. These algebras are
obtained as a finite product of totally ordered Wajsberg algebras.}$\Box
\medskip $

We denote these algebras with\textit{\ }$\left( \mathcal{W}_{i}^{n},\nabla
_{i},\rceil _{i},\mathbf{1}_{i},\leq _{i}^{n}\right) $, where $\leq _{i}^{n}$%
is the corresponding order relation on $\mathcal{W}_{i}^{n}$, $i\in
\{1,2,...,\pi _{n}\}$.\medskip \textit{\ }

\textbf{Remark 4.9.} We will denote with $\left( \mathcal{W}_{ij}^{n},\nabla
_{ij},\rceil _{ij},\mathbf{1}_{i},\leq _{ij}\right) $ the Wajsberg algebras
isomorphic to $\mathcal{W}_{i}^{n}$, as ordered sets, where $\leq _{ij}^{n}$%
is the corresponding order relation on $\mathcal{W}_{ij}^{n}$. Let $%
f_{ij}^{n}:\mathcal{W}_{i}^{n}\rightarrow \mathcal{W}_{ij}^{n}$ be such an
isomorphism of ordered sets. The Wajsberg structure on the algebra $\mathcal{%
W}_{ij}^{n}$ is given as follows. Let $x,y\in $ $\mathcal{W}_{ij}^{n}$ and $%
a,b\in \mathcal{W}_{i}^{n}$ such that $f_{ij}^{n}\left( a\right) =x$ and $%
f_{ij}^{n}\left( b\right) =y$. We define%
\begin{equation*}
x\nabla _{ij}y=f_{ij}^{n}\left( a\right) \nabla _{ij}f_{ij}^{n}\left(
b\right) \overset{def}{=}f_{ij}^{n}\left( a\nabla _{i}b\right) \text{.}
\end{equation*}%
We remark that this is the only way to define a Wajsberg algebra structure
on $\mathcal{W}_{ij}^{n}$ such that the induced order relation on this
algebra is $\nabla _{ij}$. Algebras $\mathcal{W}_{i}^{n}$ and $\mathcal{W}%
_{ij}^{n}$ are isomorphic as ordered sets and are not always isomorphic as
Wajsberg algebras. 

From the above we obtain an algorithm to find all finite Wajsberg algebras
of order $n$.\medskip 

\textbf{Example 4.10. }

i) Let $n=6$. We have that $n=2\cdot 3=3\cdot 2$, therefore $\pi _{6}=1$.

ii) Let $n=8$. We have that $n=2\cdot 4=4\cdot 2=2\cdot 2\cdot 2$.
Therefore, $\pi _{8}=2$.

iii) Let $n=12$. We have that $n=2\cdot 2\cdot 3=2\cdot 3\cdot 2=$\newline
$=3\cdot 2\cdot 2=4\cdot 3=3\cdot 4=2\cdot 6=6\cdot 2$. Then, $\pi _{12}=3$%
.\medskip

\textbf{Example 4.11}. There is only one type of partially ordered Wajsberg
algebra with $4$ elements, up to an isomorphism. Indeed, let $\left(
W_{1}=\{0,1\},\circ ,\overline{},1\right) $ and $\left( W_{2}=\{0,e\},\cdot
,^{\prime },e\right) $ be two finite totally ordered Wajsberg algebras. We
consider $W_{1}\times W_{2}=\{\left( 0,0\right) ,\left( 0,e\right) ,\left(
1,0\right) ,\left( 1,e\right) \}=$\newline
$=\{O,A,B,E\}$. On $W_{1}\times W_{2}$ we obtain a Wajsberg algebra
structure by defining the multiplication as in relation $\left( 4.2\right) $%
. We give this multiplication in the following table:

\begin{equation*}
\begin{tabular}{l|llll}
$\nabla $ & $O$ & $A$ & $B$ & $E$ \\ \hline
$O$ & $E$ & $E$ & $E$ & $E$ \\ 
$A$ & $B$ & $E$ & $B$ & $E$ \\ 
$B$ & $A$ & $A$ & $E$ & $E$ \\ 
$E$ & $O$ & $A$ & $B$ & $E$%
\end{tabular}%
\text{.}
\end{equation*}

\textbf{Example 4.12.} There is only one type of partially ordered Wajsberg
algebra with $6$ elements, up to an isomorphism. Indeed, $\pi _{6}=1$. Let $%
\left( W_{1}=\{0,1\},\circ ,\overline{},1\right) $ and $\left(
W_{2}=\{0,b,e\},\cdot ,^{\prime },e\right) $ be two finite totally ordered
Wajsberg algebras. Using relation $\left( 4.1\right) $, on $W_{2}$ we have
that $b^{\prime }=b$. We consider $W_{1}\times W_{2}=\{\left( 0,0\right)
,\left( 0,b\right) ,\left( 0,e\right) ,\left( 1,0\right) ,\left( 1,b\right)
,\left( 1,e\right) \}=$\newline
$=\{O,A,B,C,D,E\}$. On $W_{1}\times W_{2}$ we obtain a Wajsberg algebra
structure by defining the multiplication as in relation $\left( 4.2\right) $%
. We give this multiplication in the following table:

\begin{equation}
\begin{tabular}{l|llllll}
$\nabla _{11}^{6}$ & $O$ & $A$ & $B$ & $C$ & $D$ & $E$ \\ \hline
$O$ & $E$ & $E$ & $E$ & $E$ & $E$ & $E$ \\ 
$A$ & $D$ & $E$ & $E$ & $D$ & $E$ & $E$ \\ 
$B$ & $C$ & $D$ & $E$ & $C$ & $D$ & $E$ \\ 
$C$ & $B$ & $B$ & $B$ & $E$ & $E$ & $E$ \\ 
$D$ & $A$ & $B$ & $B$ & $D$ & $E$ & $E$ \\ 
$E$ & $O$ & $A$ & $B$ & $C$ & $D$ & $E$%
\end{tabular}%
\ .  \tag{4.4.}
\end{equation}

We remark that $A\leq B,A\leq D,C\leq D$, and the other elements can't be
compared in the algebra $\mathcal{W}_{1}^{6}=\mathcal{W}_{11}^{6}=\left(
W_{1}\times W_{2},\nabla _{11}^{6}\right) $. We denote this order relation
with $\leq _{11}^{6}$.

If we consider the isomorphism 
\begin{eqnarray*}
f_{12}^{6} &:&\left( W_{1}\times W_{2},\nabla _{11}^{6}\right) \rightarrow
\left( W_{1}\times W_{2},\nabla _{12}^{6}\right) , \\
f_{12}^{6}\left( A\right)  &\text{=}&A,f_{12}^{6}\left( B\right) \text{=}%
C,f_{12}^{6}\left( C\right) \text{=}B,f_{12}^{6}\left( D\right) \text{=}%
D,f_{12}^{6}\left( O\right) \text{=}O,f_{12}^{6}\left( E\right) \text{=}E%
\text{,}
\end{eqnarray*}%
we obtain on $W_{1}\times W_{2}$ a new Wajsberg algebra structure, with the
multiplication $\nabla _{12}^{6}$ given in the below table:

\begin{equation}
\begin{tabular}{l|llllll}
$\nabla _{12}^{6}$ & $O$ & $A$ & $B$ & $C$ & $D$ & $E$ \\ \hline
$O$ & $E$ & $E$ & $E$ & $E$ & $E$ & $E$ \\ 
$A$ & $D$ & $E$ & $D$ & $E$ & $E$ & $E$ \\ 
$B$ & $C$ & $C$ & $E$ & $C$ & $E$ & $E$ \\ 
$C$ & $B$ & $D$ & $B$ & $E$ & $D$ & $E$ \\ 
$D$ & $A$ & $C$ & $D$ & $C$ & $E$ & $E$ \\ 
$E$ & $O$ & $A$ & $B$ & $C$ & $D$ & $E$%
\end{tabular}
\tag{4.5.}
\end{equation}

These two algebras, $\mathcal{W}_{11}^{6}=\left( W_{1}\times W_{2},\nabla
_{11}^{6}\right) $ and $\mathcal{W}_{12}^{6}=\left( W_{1}\times W_{2},\nabla
_{12}^{6}\right) $, are isomorphic. We remark that $A\leq C,A\leq D,B\leq D$
and the other elements can't be compared in the algebra $\mathcal{W}_{12}^{6}
$. We denote this order relation with $\leq _{12}^{6}$. This \ algebra is
the Wajsberg algebra given in Example 3.18. The isomorphism $f_{12}^{6}$ is
in the same time isomorphism of ordered sets and isomorphism of Wajsberg
algebras.

If we consider the isomorphism 
\begin{eqnarray*}
f_{13}^{6} &:&\left( W_{1}\times W_{2},\nabla _{11}^{6}\right) \rightarrow
\left( W_{1}\times W_{2},\nabla _{13}^{6}\right) , \\
f_{13}^{6}\left( A\right)  &\text{=}&B,f_{13}^{6}\left( B\right) \text{=}%
D,f_{13}^{6}\left( C\right) \text{=}C,f_{13}^{6}\left( D\right) \text{=}%
A,f_{13}^{6}\left( O\right) \text{=}O,f_{13}^{6}\left( E\right) \text{=}E,
\end{eqnarray*}%
we obtain on $W_{1}\times W_{2}$ a new Wajsberg algebra structure, with the
multiplication $\nabla _{13}^{6}$ given in the below table: 
\begin{equation}
\begin{tabular}{l|llllll}
$\nabla _{13}^{6}$ & $O$ & $A$ & $B$ & $C$ & $D$ & $E$ \\ \hline
$O$ & $E$ & $E$ & $E$ & $E$ & $E$ & $E$ \\ 
$A$ & $B$ & $E$ & $D$ & $A$ & $D$ & $E$ \\ 
$B$ & $A$ & $E$ & $E$ & $A$ & $E$ & $E$ \\ 
$C$ & $D$ & $E$ & $D$ & $E$ & $D$ & $E$ \\ 
$D$ & $C$ & $A$ & $A$ & $C$ & $E$ & $E$ \\ 
$E$ & $O$ & $A$ & $B$ & $C$ & $D$ & $E$%
\end{tabular}%
.  \tag{4.6.}
\end{equation}%
In $\mathcal{W}_{13}^{6}=\left( W_{1}\times W_{2},\nabla _{13}^{6}\right) $,
we have $B\leq A,C\leq A,B\leq D$ and the other elements can't be compared.
We denote this order relation with $\leq _{13}^{6}$.The isomorphism $%
f_{13}^{6}$ is only isomorphism of ordered sets and is not isomorphism of
Wajsberg algebras.

If we consider $W_{2}\times W_{1}=\{\left( 0,0\right) ,\left( 0,1\right)
,\left( b,0\right) ,\left( b,1\right) ,\left( e,0\right) ,\left( e,1\right)
\}=$\newline
$=\{O,A,B,C,D,E\}$, on $W_{2}\times W_{1}$ we obtain a Wajsberg algebra
structure by defining the multiplication as in relation $\left( 4.2\right) $%
. We give this multiplication in the following table:

\begin{equation}
\begin{tabular}{l|llllll}
$\nabla _{14}^{6}$ & $O$ & $A$ & $B$ & $C$ & $D$ & $E$ \\ \hline
$O$ & $E$ & $E$ & $E$ & $E$ & $E$ & $E$ \\ 
$A$ & $D$ & $E$ & $D$ & $E$ & $D$ & $E$ \\ 
$B$ & $C$ & $C$ & $E$ & $E$ & $E$ & $E$ \\ 
$C$ & $B$ & $C$ & $D$ & $E$ & $D$ & $E$ \\ 
$D$ & $A$ & $A$ & $C$ & $C$ & $E$ & $E$ \\ 
$E$ & $O$ & $A$ & $B$ & $C$ & $D$ & $E$%
\end{tabular}%
.  \tag{4.7.}
\end{equation}

The algebras $\left( W_{1}\times W_{2},\nabla _{11}^{6}\right) $ and $\left(
W_{2}\times W_{1},\nabla _{14}^{6}\right) $ are also isomorphic, by taking
the map 
\begin{eqnarray*}
f_{14}^{6} &:&\left( W_{1}\times W_{2},\nabla _{11}^{6}\right) \rightarrow
\left( W_{2}\times W_{1},\nabla _{14}^{6}\right) , \\
f_{14}^{6}\left( A\right)  &\text{=}&B,f_{14}^{6}\left( B\right) \text{=}%
D,f_{14}^{6}\left( C\right) \text{=}A,f_{14}^{6}\left( D\right) \text{=}%
C,f_{14}^{6}\left( O\right) \text{=}O,f_{14}^{6}\left( E\right) \text{=}E.
\end{eqnarray*}%
In $\mathcal{W}_{14}^{6}=\left( W_{2}\times W_{1},\nabla _{14}^{6}\right) $,
we have $A\leq C,B\leq C,B\leq D~$and the other elements can't be compared.
We denote this order relation with $\leq _{14}^{6}$.The isomorphism $%
f_{14}^{6}$ is in the same time isomorphism of ordered sets and isomorphism
of Wajsberg algebras.\medskip 

\textbf{Example 4.13}. There is only two types of partially ordered Wajsberg
algebra with $8$ elements, up to an isomorphism. Indeed, $\pi _{8}=2$. Let 
\newline
$\left( W_{1}=\{0,a,b,e\},\circ ,\overline{},1\right) $ and $\left(
W_{2}=\{0,1\},\cdot ,^{\prime },e\right) $ be two finite totally ordered
Wajsberg algebras. Using relation $\left( 4.1\right) $, on $W_{1}$ we have
that $\overline{b}=a$ and $\overline{a}=b$. We consider $W_{1}\times
W_{2}=\{\left( 0,0\right) ,\left( 0,1\right) ,\left( a,0\right) ,\left(
a,1\right) ,\left( b,0\right) ,\left( b,1\right) ,\left( e,0\right) ,\left(
e,1\right) \}=$\newline
$=\{O,X,Y,Z,T,U,V,E\}$. On $W_{1}\times W_{2}$ we obtain a Wajsberg algebra
structure by defining the multiplication as in relation $\left( 4.2\right) $%
, namely $\mathcal{W}_{11}^{8}=\left( W_{1}\times W_{2},\nabla
_{11}^{8}\right) $. The multiplication $\nabla _{11}^{8}$ is given in the
following table:%
\begin{equation}
\begin{tabular}{l|llllllll}
$\nabla _{11}^{8}$ & $O$ & $X$ & $Y$ & $Z$ & $T$ & $U$ & $V$ & $E$ \\ \hline
$O$ & $E$ & $E$ & $E$ & $E$ & $E$ & $E$ & $E$ & $E$ \\ 
$X$ & $V$ & $E$ & $V$ & $E$ & $V$ & $E$ & $V$ & $E$ \\ 
$Y$ & $U$ & $U$ & $E$ & $E$ & $E$ & $E$ & $E$ & $E$ \\ 
$Z$ & $T$ & $U$ & $V$ & $E$ & $V$ & $E$ & $V$ & $E$ \\ 
$T$ & $Z$ & $Z$ & $U$ & $U$ & $E$ & $E$ & $E$ & $E$ \\ 
$U$ & $Y$ & $Z$ & $T$ & $U$ & $T$ & $E$ & $V$ & $E$ \\ 
$V$ & $X$ & $X$ & $Z$ & $Z$ & $U$ & $U$ & $E$ & $E$ \\ 
$E$ & $O$ & $X$ & $Y$ & $Z$ & $T$ & $U$ & $V$ & $E$%
\end{tabular}
\tag{4.8.}
\end{equation}

In $\mathcal{W}_{11}^{8}$ we have that $O\leq X\leq Z\leq U\leq E$, $O\leq
Y\leq T\leq V\leq E$, $O\leq Y\leq Z\leq U\leq E,O\leq Y\leq T\leq U\leq E$
and the other elements can't be compared in this algebra. We denote this
order relation with $\leq _{11}^{8}$.

Now, we consider $W_{2}\times W_{1}=\{\left( 0,0\right) ,\left( 0,a\right)
,\left( 0,b\right) ,\left( 0,e\right) ,\left( 1,0\right) ,\left( 1,a\right)
,\left( 1,b\right) ,\left( 1,e\right) \}=$\newline
$=\{O,X,Y,Z,T,U,V,E\}$. On $W_{2}\times W_{1}$ we obtain a Wajsberg algebra
structure by defining the multiplication as in relation $\left( 4.2\right) $%
, namely $\mathcal{W}_{12}^{8}=\left( W_{2}\times W_{1},\nabla
_{12}^{8}\right) $. The multiplication $\nabla _{12}^{8}$ is given in the
following table:%
\begin{equation}
\begin{tabular}{l|llllllll}
$\nabla _{12}^{8}$ & $O$ & $X$ & $Y$ & $Z$ & $T$ & $U$ & $V$ & $E$ \\ \hline
$O$ & $E$ & $E$ & $E$ & $E$ & $E$ & $E$ & $E$ & $E$ \\ 
$X$ & $V$ & $E$ & $E$ & $E$ & $V$ & $E$ & $E$ & $E$ \\ 
$Y$ & $U$ & $V$ & $E$ & $E$ & $U$ & $V$ & $E$ & $E$ \\ 
$Z$ & $T$ & $U$ & $V$ & $E$ & $T$ & $U$ & $V$ & $E$ \\ 
$T$ & $Z$ & $X$ & $Z$ & $Z$ & $E$ & $E$ & $E$ & $E$ \\ 
$U$ & $Y$ & $Z$ & $Z$ & $Z$ & $V$ & $E$ & $E$ & $E$ \\ 
$V$ & $X$ & $Y$ & $Z$ & $Z$ & $U$ & $V$ & $E$ & $E$ \\ 
$E$ & $O$ & $X$ & $Y$ & $Z$ & $T$ & $U$ & $V$ & $E$%
\end{tabular}
\tag{4.8.}
\end{equation}%
We have that $O\leq X\leq Y\leq Z\leq E$, $O\leq X\leq Y\leq V\leq E$,%
\newline
$O\leq X\leq U\leq V\leq E$, $O\leq T\leq U\leq V\leq E$. These two
structures, $\mathcal{W}_{11}^{8}$ and $\mathcal{W}_{12}^{8}$, are
isomorphic. The morphism is 
\begin{eqnarray*}
f_{12}^{8} &:&\left( W_{1}\times W_{2},\nabla _{11}^{8}\right) \rightarrow
\left( W_{2}\times W_{1},\nabla _{12}^{8}\right) , \\
f_{12}^{8}\left( X\right)  &\text{=}&T,f_{12}^{8}\left( Y\right) \text{=}%
X,f_{12}^{8}\left( Z\right) \text{=}U,f_{12}^{8}\left( T\right) \text{=}Y, \\
f_{12}^{8}\left( U\right)  &\text{=}&V,f_{12}^{8}\left( V\right) \text{=}%
Z,f_{12}^{8}\left( O\right) \text{=}O,f_{12}^{8}\left( E\right) \text{=}E.
\end{eqnarray*}%
The isomorphism $f_{12}^{8}$ is in the same time isomorphism of ordered sets
and isomorphism of Wajsberg algebras.

If we take $W_{2}\times W_{2}\times W_{2}=\{\left( 0,0,0\right) ,\left(
0,0,1\right) ,\left( 0,1,0\right) ,\left( 0,1,1\right) ,\left( 1,0,0\right) ,
$\newline
$\left( 1,0,1\right) ,\left( 1,1,0\right) ,\left( 1,1,1\right)
\}=\{O,X,Y,Z,T,U,V,E\}$, on $W_{2}\times W_{2}\times W_{2}$ we obtain a
Wajsberg algebra structure by defining the multiplication as in relation $%
\left( 4.2\right) $, namely $\mathcal{W}_{21}^{8}=\left( W_{2}\times
W_{2}\times W_{2},\nabla _{21}^{8}\right) $. The multiplication $\nabla
_{21}^{8}$ is given in the following table:%
\begin{equation}
\begin{tabular}{l|llllllll}
$\nabla _{21}^{8}$ & $O$ & $X$ & $Y$ & $Z$ & $T$ & $U$ & $V$ & $E$ \\ \hline
$O$ & $E$ & $E$ & $E$ & $E$ & $E$ & $E$ & $E$ & $E$ \\ 
$X$ & $V$ & $E$ & $V$ & $E$ & $V$ & $E$ & $V$ & $E$ \\ 
$Y$ & $U$ & $U$ & $E$ & $E$ & $U$ & $U$ & $E$ & $E$ \\ 
$Z$ & $T$ & $U$ & $V$ & $E$ & $T$ & $U$ & $V$ & $E$ \\ 
$T$ & $Z$ & $Z$ & $Z$ & $Z$ & $E$ & $E$ & $E$ & $E$ \\ 
$U$ & $Y$ & $Z$ & $Y$ & $Z$ & $V$ & $E$ & $V$ & $E$ \\ 
$V$ & $X$ & $X$ & $Z$ & $Z$ & $U$ & $U$ & $E$ & $E$ \\ 
$E$ & $O$ & $X$ & $Y$ & $Z$ & $T$ & $U$ & $V$ & $E$%
\end{tabular}
\tag{4.9.}
\end{equation}%
We have that $X\leq Z,X\leq U,Y\leq Z,T\leq V,Y\leq V,T\leq U$. These two
structures, $\mathcal{W}_{11}^{8}$ and $\mathcal{W}_{21}^{8}$, are not
isomorphic as ordered sets, therefore are not isomorphic as Wajsberg
algebras.\medskip 

\textbf{Remark 4.14.} 1) There is only one type of partially ordered
Wajsberg algebra with $9$ elements, up to an isomorphism as ordered sets.

2) There is only one type of partially ordered Wajsberg algebra with $10$
elements, up to an isomorphism, up to an isomorphism as ordered sets..

3) There is only three types of partially ordered Wajsberg algebra with $12$
elements, up to an isomorphism, up to an isomorphism as ordered sets..

\begin{equation*}
\end{equation*}

\textbf{5. Special algebras arising from block codes }%
\begin{equation*}
\end{equation*}

In [FL; 15], was developed an algorithm which provide conditions to attach a
BCK algebra to a given block code. In the following, we will use those ideas
to obtain a similar algorithm in the case of MV-algebras and Wajsberg
algebras. The difference is that, in the first case, the BCK algebra arising
from a block code is a non-commutative, a non-implicative, but a positive
implicative BCK algebra ([FL; 17]). In the second case, of MV-algebras and
Wajsberg algebras, we must obtain from a given block-code a BCK commutative
bounded algebra. Our idea is to obtain a Wajsberg finite algebra associated
to a given block code and from here the desired MV-algebra and the desired
bounded commutative BCK-algebra.

Let $V$ be a binary block-code with $n+1$ codewords of length $n+1$,
lexicographically ordered, $V=\{w_{0},w_{1},...,w_{n}\}$ and $%
w_{x}=x_{1}x_{2}...x_{n}x_{n+1}\in V,$ $w_{y}=y_{1}y_{2}...y_{n}y_{n+1}\in V$
be two codewords. On $V$ we can define the following partial order relation: 
\begin{equation}
w_{x}\preceq w_{y}\text{ if and only if }y_{i}\leq x_{i},i\in
\{0,1,2,...,n\}.  \tag{5.1.}
\end{equation}

We consider the matrix $M_{V}=\left( m_{i,j}\right) _{i,j\in
\{1,2,...,n+1\}}\in \mathcal{M}_{n+1}(\{0,1\})$ with the rows consisting of
the codewords of $V.$ This matrix is called \textit{the matrix associated to
the block-code} $V.\medskip $

\textbf{Theorem 5.1.} \textit{With the above notations, if \ the matrix }$%
M_{V}$ \textit{has the first line and the last column of the form} $\underset%
{n+1-\text{time}}{\underbrace{11...1}}$, \textit{the last line of the form} $%
\underset{n-\text{time}}{\underbrace{00...00}1}$\textit{, the first column
of the form} $\underset{n-\text{time}}{1\underbrace{00...00}}$, $m_{ii}=1$, 
\textit{for all} $i\in \{1,2,...,n+1\}$, \textit{and if the order relation
given by }$\left( 5.1\right) $ \textit{coincide with one of the order
relations }$\leq _{ij}^{n+1}$\textit{given by the Remark 4.9, then} \textit{%
there are a set} $A$ \textit{with} $n+1$ \textit{elements, a Wajsberg
algebra }$\mathcal{W}_{ij}^{n+1}$ \textit{and a W-function} $f:A\rightarrow 
\mathcal{W}_{ij}^{n+1}$ \textit{such that} $f$ \textit{determines} $%
V.\medskip $

\textbf{Proof. }\ We consider on $V$ the lexicographic order, denoted by $%
\leq _{lex}$. \ It results that $(V,\leq _{lex})$ is a totally ordered set.
Let $V=\{w_{0},w_{1},...,w_{n}\},$ with $w_{n}\geq _{lex}w_{n-1}\geq
_{lex}...\geq _{lex}w_{0}$. We denote $w_{n}=\underset{n+1-\text{time}}{%
\underbrace{11...1}}$ and $w_{0}=\underset{n-\text{time}}{\underbrace{00...0}%
1}$. We remark that $w_{0}=\mathbf{1}$ is the "one" element and $%
w_{n}=\theta $ is the "zero" element in $\left( V,\preceq \right) $,
considered with order relation $\left( 5.1\right) $. If this order relation
coincides with one of the order relation $\leq _{ij}^{n+1}$, given in Remark
4.9, it results that on $\left( V,\preceq \right) $ we can obtain a Wajsberg
algebra structure, which is isomorphic to a Wajsberg algebra $\mathcal{W}%
_{ij}^{n+1}$, with attached order relation $\leq _{ij}^{n+1}$. If we
consider $A=V$ and the identity map $f:A\rightarrow V,f\left( w\right) =w$
as a W-function, the decomposition of \ $f$ provides a family of maps $V_{%
\mathcal{W}_{ij}^{n+1}}=\{f_{r}:A\rightarrow \{0,1\}~/~$ \ $f_{r}\left(
x\right) =1$ if and only if \ $r\ast f\left( x\right) =w_{0},\forall x\in
A,r\in X\}.$ This family is the binary block-code $V$ relative to the order
relation $\preceq .\medskip \Box $

The above Theorem extends to W-algebras results obtained in Theorem 3.2 from
[FL; 15]\medskip .

\textbf{Proposition 5.2. }\ \textit{Let} \ $A=\left( a_{i,j}\right) 
_{\substack{ i\in \{1,2,...,n\}  \\ j\in \{1,2,...,m\}}}\in \mathcal{M}%
_{n,m}(\{0,1\})$ \textit{be a matrix. Starting from this matrix, we can find
a \ matrix} $B=\left( b_{i,j}\right) _{i,j\in \{1,2,...,q\}}\in \mathcal{M}%
_{q}(\{0,1\})$, $q\geq max\{m,n\}$, \textit{such that} $B$ \textit{has the
first line and the last column of the form} $\underset{q-\text{time}}{%
\underbrace{11...1}}$, \textit{the last line of the form} $\underset{q-1-%
\text{time}}{\underbrace{00...00}1}$, \textit{the first column of the form} $%
\underset{q-1-\text{time}}{1\underbrace{00...00}}~$\textit{with }$%
b_{ii}=1,\forall i\in \{1,2,...,q\}$ \textit{and} $A$ \textit{becomes a
submatrix of the matrix} $B.\medskip $

\textbf{Proof. }Obviously$.\Box \medskip $

\textbf{Theorem 5.3.} \textit{With the above notations, let } $V$ be\textit{%
\ a binary block-code with} $n$ \textit{codewords of length} $m,n\neq m$. 
\textit{Let} $q$ \textit{be} \textit{a natural number, }$q\geq max\{m,n\}$, 
\textit{and} $B=\left( b_{i,j}\right) _{i,j\in \{1,2,...,q\}}\in \mathcal{M}%
_{q}(\{0,1\})$ \textit{be a matrix such that the matrix} $M_{V}$ \textit{is a%
} \textit{submatrix of the matrix} $B$. \textit{If} $B$ \textit{is the
matrix attached to a code} $C$ \textit{which satisfies the conditions from
Theorem 5.1, therefore there is} \textit{a set} $A$ \textit{with} $m$ 
\textit{elements, a Wajsberg algebra }$\mathcal{W}_{ij}^{n+1}$\textit{\ and
a W-function} $f:A\rightarrow \mathcal{W}_{ij}^{n+1}$ \textit{such that the
obtained block-code }$C_{nm}$ \textit{contains the block-code} $V$ \textit{%
as a subset.\medskip }

\textbf{Proof.} Let $V$ be a binary block-code, $V=\{w_{1,}w_{2},...,w_{n}\}$%
, with codewords of length $m$. Let $M\in \mathcal{M}_{n,m}(\{0,1\})$ be its
associated matrix. Using Proposition 5.2, we can extend the matrix $M$ to a
square matrix $M^{\prime }\in \mathcal{M}_{q}(\{0,1\})$ and we can apply
Theorem 5.1 for the matrix $M^{\prime }$.$\ $Assuming that the initial
columns of the matrix $M$ \ have in the new matrix $M^{\prime }$ positions $%
i_{j_{1}},i_{j_{2}},...,i_{j_{m}}\in \{1,2,...,q\},$ let $%
A=\{x_{j_{1}},x_{j_{2}},...,x_{j_{m}}\}\subseteq C.$ The W-function $%
f:A\rightarrow C,f\left( x_{j_{i}}\right) =$ $x_{j_{i}},$ $i\in
\{1,2,...,m\} $, determines the binary block-code $C_{nm}$ such that $%
V\subseteq C_{nm}.\Box \medskip $

The above Theorem extends to W-algebras results obtained in Theorem 3.9 from
[FL; 15]\medskip .\medskip

\textbf{Remark 5.4.} In the above theorem, the associated Wajsberg algebra
is not unique, as we can see in the Example 6.10.\medskip

\textbf{Theorem 5.5.} \textit{If a block-code} $V$ \textit{satisfies the
conditions from Theorem 5.1, there is an MV-algebra} $\mathcal{X}$ \textit{%
and a commutative bounded BCK-algebra} $\mathcal{Y}$ \textit{associated to
this code.\medskip }

\textbf{Proof.} We use Theorem 5.1, Remark 2.4 and Remark 2.6.%
\begin{equation*}
\end{equation*}

\begin{equation*}
\end{equation*}

\textbf{6. Examples}%
\begin{equation*}
\end{equation*}

\textbf{Example 6.1.} We consider the code \newline
$V_{1}=\{111111,011011,001001,000111,000011,000001\}=$\newline
$=\{\theta ,a,b,c,d,e\}$, with elements lexicographically ordered. The
associated matrix is%
\begin{equation*}
\begin{tabular}{l|llllll}
& $\theta $ & $a$ & $b$ & $c$ & $d$ & $e$ \\ \hline
$\theta $ & $\mathbf{1}$ & $\mathbf{1}$ & $\mathbf{1}$ & $\mathbf{1}$ & $%
\mathbf{1}$ & $\mathbf{1}$ \\ 
$a$ & $0$ & $\mathbf{1}$ & $\mathbf{1}$ & $0$ & $\mathbf{1}$ & $\mathbf{1}$
\\ 
$b$ & $0$ & $0$ & $\mathbf{1}$ & $0$ & $0$ & $\mathbf{1}$ \\ 
$c$ & $0$ & $0$ & $0$ & $\mathbf{1}$ & $\mathbf{1}$ & $\mathbf{1}$ \\ 
$d$ & $0$ & $0$ & $0$ & $0$ & $\mathbf{1}$ & $\mathbf{1}$ \\ 
$e$ & $0$ & $0$ & $0$ & $0$ & $0$ & $\mathbf{1}$%
\end{tabular}%
.
\end{equation*}%
Using order \ given by the relation $\left( 5.1\right) $, we have that $%
a\preceq b,c\preceq d,a\preceq d$ and the other elements can't be compared.
Therefore this order relation is $\leq _{11}^{6}$ and the associated
Wajsberg algebra is algebra $\mathcal{W}_{11}^{6}$, given by relation $%
\left( 4.4\right) $.\medskip

\textbf{Example 6.2. }We consider the code \newline
$V_{2}=\{111111,010111,001011,000101,000011,000001\}=$\newline
$=\{\theta ,a,b,c,d,e\}$, with elements lexicographically ordered. The
associated matrix is%
\begin{equation*}
~%
\begin{tabular}{l|llllll}
& $\theta $ & $a$ & $b$ & $c$ & $d$ & $e$ \\ \hline
$\theta $ & $\mathbf{1}$ & $\mathbf{1}$ & $\mathbf{1}$ & $\mathbf{1}$ & $%
\mathbf{1}$ & $\mathbf{1}$ \\ 
$a$ & $0$ & $\mathbf{1}$ & $0$ & $\mathbf{1}$ & $\mathbf{1}$ & $\mathbf{1}$
\\ 
$b$ & $0$ & $0$ & $\mathbf{1}$ & $0$ & $\mathbf{1}$ & $\mathbf{1}$ \\ 
$c$ & $0$ & $0$ & $0$ & $\mathbf{1}$ & $0$ & $\mathbf{1}$ \\ 
$d$ & $0$ & $0$ & $0$ & $0$ & $\mathbf{1}$ & $\mathbf{1}$ \\ 
$e$ & $0$ & $0$ & $0$ & $0$ & $0$ & $\mathbf{1}$%
\end{tabular}%
.
\end{equation*}%
Using order \ given by the relation $\left( 5.1\right) ,$ we have that $%
a\preceq c,a\preceq d,b\preceq d$ and the other elements can't be compared.
Therefore this order relation is $\leq _{12}^{6}$ and the associated
Wajsberg algebra is $\mathcal{W}_{12}^{6}$, given by relation $\left(
4.5\right) $. Algebra $\mathcal{W}_{12}^{6}$ is the Wajsberg algebra given
in the Example 3.18. From here, we can see the MV-algebra and the BCK
commutative bounded algebra associated to the block-code $V_{2}$.\medskip

\textbf{Example 6.3.} We consider the code \newline
$V_{3}=\{111111,010001,011011,010101,000011,000001\}=$\newline
$=\{\theta ,a,b,c,d,e\}$, with elements lexicographically ordered. The
associated matrix is%
\begin{equation*}
~%
\begin{tabular}{l|llllll}
& $\theta $ & $a$ & $b$ & $c$ & $d$ & $e$ \\ \hline
$\theta $ & $\mathbf{1}$ & $\mathbf{1}$ & $\mathbf{1}$ & $\mathbf{1}$ & $%
\mathbf{1}$ & $\mathbf{1}$ \\ 
$a$ & $0$ & $\mathbf{1}$ & $0$ & $0$ & $0$ & $\mathbf{1}$ \\ 
$b$ & $0$ & $\mathbf{1}$ & $\mathbf{1}$ & $0$ & $\mathbf{1}$ & $\mathbf{1}$
\\ 
$c$ & $0$ & $\mathbf{1}$ & $0$ & $\mathbf{1}$ & $0$ & $\mathbf{1}$ \\ 
$d$ & $0$ & $0$ & $0$ & $0$ & $\mathbf{1}$ & $\mathbf{1}$ \\ 
$e$ & $0$ & $0$ & $0$ & $0$ & $0$ & $\mathbf{1}$%
\end{tabular}%
.
\end{equation*}%
Using order \ given by the relation $\left( 5.1\right) ,$ we have that $%
b\preceq a,c\preceq a,b\preceq d$ and the other elements can't be compared.
Therefore this order relation is $\leq _{13}^{6}$ and the associated
Wajsberg algebra is $\mathcal{W}_{13}^{6}$, given by relation $\left(
4.6\right) $.

The skeleton of this algebra is

\begin{equation}
\begin{tabular}{l|llllll}
& $\theta $ & $a$ & $b$ & $c$ & $d$ & $e$ \\ \hline
$\theta $ & $\mathbf{\blacksquare }$ & $\mathbf{\blacksquare }$ & $\mathbf{%
\blacksquare }$ & $\mathbf{\blacksquare }$ & $\mathbf{\blacksquare }$ & $%
\mathbf{\blacksquare }$ \\ 
$a$ &  & $\mathbf{\blacksquare }$ &  &  &  & $\mathbf{\blacksquare }$ \\ 
$b$ &  & $\mathbf{\blacksquare }$ & $\mathbf{\blacksquare }$ &  & $\mathbf{%
\blacksquare }$ & $\mathbf{\blacksquare }$ \\ 
$c$ &  & $\mathbf{\blacksquare }$ &  & $\mathbf{\blacksquare }$ &  & $%
\mathbf{\blacksquare }$ \\ 
$d$ &  &  &  &  & $\mathbf{\blacksquare }$ & $\mathbf{\blacksquare }$ \\ 
$e$ &  &  &  &  &  & $\mathbf{\blacksquare }$%
\end{tabular}%
.  \tag{6.1.}
\end{equation}

\textbf{Example 6.4.} We consider the code \newline
$V_{4}=\{111111,010101,001111,000101,000011,000001\}=$\newline
$=\{\theta ,a,b,c,d,e\}$, with elements lexicographically ordered. Using
order \ given by the relation $\left( 5.1\right) ,$ we have that $a\preceq
c,b\preceq c,b\preceq d$ and the other elements can't be compared. Therefore
this order relation is $\leq _{14}^{6}$ and the associated Wajsberg algebra
is $\mathcal{W}_{14}^{6}$, given by relation $\left( 4.7\right) $.

The skeleton of this algebra is

\begin{equation}
\begin{tabular}{l|llllll}
& $\theta $ & $a$ & $b$ & $c$ & $d$ & $e$ \\ \hline
$\theta $ & $\mathbf{\blacksquare }$ & $\mathbf{\blacksquare }$ & $\mathbf{%
\blacksquare }$ & $\mathbf{\blacksquare }$ & $\mathbf{\blacksquare }$ & $%
\mathbf{\blacksquare }$ \\ 
$a$ &  & $\mathbf{\blacksquare }$ &  & $\mathbf{\blacksquare }$ &  & $%
\mathbf{\blacksquare }$ \\ 
$b$ &  &  & $\mathbf{\blacksquare }$ & $\mathbf{\blacksquare }$ & $\mathbf{%
\blacksquare }$ & $\mathbf{\blacksquare }$ \\ 
$c$ &  &  &  & $\mathbf{\blacksquare }$ &  & $\mathbf{\blacksquare }$ \\ 
$d$ &  &  &  &  & $\mathbf{\blacksquare }$ & $\mathbf{\blacksquare }$ \\ 
$e$ &  &  &  &  &  & $\mathbf{\blacksquare }$%
\end{tabular}%
.  \tag{6.2.}
\end{equation}

\textbf{Example 6.5. \ }We consider the code $V_{5}=\{011,101,010,001,000\}$%
. Using Theorem 5.3, we have the matrix $B\,\ $ 
\begin{equation*}
B=~%
\begin{tabular}{|llllll}
$\theta $ & $a$ & $b$ & $c$ & $d$ & $e$ \\ \hline
$1$ & $1$ & $\mathbf{1}$ & $\mathbf{1}$ & $\mathbf{1}$ & $1$ \\ 
$0$ & $1$ & $\mathbf{0}$ & $\mathbf{1}$ & $\mathbf{1}$ & $1$ \\ 
$0$ & $0$ & $\mathbf{1}$ & $\mathbf{0}$ & $\mathbf{1}$ & $1$ \\ 
$0$ & $0$ & $\mathbf{0}$ & $\mathbf{1}$ & $\mathbf{0}$ & $1$ \\ 
$0$ & $0$ & $\mathbf{0}$ & $\mathbf{0}$ & $\mathbf{1}$ & $1$ \\ 
$0$ & $0$ & $\mathbf{0}$ & $\mathbf{0}$ & $\mathbf{0}$ & $1$%
\end{tabular}%
\text{,}
\end{equation*}%
the set $A=\{b,c,d\}$, the Wajsberg algebra $\mathcal{W}_{12}^{6}$ and the
W-function $f:A\rightarrow \mathcal{W}_{12}^{6},f\left( b\right) =b,f\left(
c\right) =c,f\left( d\right) =d$. We have the code $C_{53}=%
\{111,011,101,010,001,000\}$. The code $V_{5}=\{011,101,010,001,000\}$ is a
subcode of the code $C_{53}$.\medskip

\textbf{Example 6.6.} We consider the code \newline
$V_{6}=\{111111,011101,001101,000111,000011,000001\}=$\newline
$=\{\theta ,a,b,c,d,e\}$, with elements lexicographically ordered. The
associated matrix is 
\begin{equation*}
\begin{tabular}{l|llllll}
& $\theta $ & $a$ & $b$ & $c$ & $d$ & $e$ \\ \hline
$\theta $ & $\mathbf{1}$ & $\mathbf{1}$ & $\mathbf{1}$ & $\mathbf{1}$ & $%
\mathbf{1}$ & $\mathbf{1}$ \\ 
$a$ & $0$ & $\mathbf{1}$ & $\mathbf{1}$ & $\mathbf{1}$ & $0$ & $\mathbf{1}$
\\ 
$b$ & $0$ & $0$ & $\mathbf{1}$ & $\mathbf{1}$ & $0$ & $\mathbf{1}$ \\ 
$c$ & $0$ & $0$ & $0$ & $\mathbf{1}$ & $\mathbf{1}$ & $\mathbf{1}$ \\ 
$d$ & $0$ & $0$ & $0$ & $0$ & $\mathbf{1}$ & $\mathbf{1}$ \\ 
$e$ & $0$ & $0$ & $0$ & $0$ & $0$ & $\mathbf{1}$%
\end{tabular}%
\end{equation*}%
and the skeleton is 
\begin{equation*}
\begin{tabular}{l|llllll}
& $\theta $ & $a$ & $b$ & $c$ & $d$ & $e$ \\ \hline
$\theta $ & $\mathbf{\blacksquare }$ & $\mathbf{\blacksquare }$ & $\mathbf{%
\blacksquare }$ & $\mathbf{\blacksquare }$ & $\mathbf{\blacksquare }$ & $%
\mathbf{\blacksquare }$ \\ 
$a$ &  & $\mathbf{\blacksquare }$ & $\mathbf{\blacksquare }$ & $\mathbf{%
\blacksquare }$ &  & $\mathbf{\blacksquare }$ \\ 
$b$ &  &  & $\mathbf{\blacksquare }$ & $\mathbf{\blacksquare }$ &  & $%
\mathbf{\blacksquare }$ \\ 
$c$ &  &  &  & $\mathbf{\blacksquare }$ & $\mathbf{\blacksquare }$ & $%
\mathbf{\blacksquare }$ \\ 
$d$ &  &  &  &  & $\mathbf{\blacksquare }$ & $\mathbf{\blacksquare }$ \\ 
$e$ &  &  &  &  &  & $\mathbf{\blacksquare }$%
\end{tabular}%
\end{equation*}

Using order \ given by the relation $\left( 5.1\right) ,$ we have that $%
a\preceq b,c\preceq d$ and the other elements can't be compared. If we
consider the binary relation given by the skeleton: $i\leq _{s}j$ if and
only if in the position $\left( i,j\right) $ we have a black square, this
relation is not an order relation. Indeed, we do not have the transitivity:
we have $a\leq _{s}b,b\leq _{s}c,c\leq _{s}d$ but $a$ and $d$ can't be
compared, therefore the sets $\left( V_{6},\leq _{s}\right) $ and $\left(
V_{6},\preceq \right) $ are not isomorphic as ordered sets, as in
Proposition 3.15. Even if the associated matrix to this code is on the form
asked in Theorem 5.1, to this code we can't associate a Wajsberg algebra,
since the order relation $\preceq $ is not on the form $\leq _{ij}^{6}$,
given by Theorem 4.8 and Remark 4.9.

We must remark that, from [FL; 15], Theorem 3.2, to the code $V_{6}$ we can
attached a BCK-algebra, namely 
\begin{equation*}
B=%
\begin{tabular}{l|llllll}
& $\theta $ & $a$ & $b$ & $c$ & $d$ & $e$ \\ \hline
$\theta $ & $\theta $ & $\theta $ & $\theta $ & $\theta $ & $\theta $ & $%
\theta $ \\ 
$a$ & $a$ & $\theta $ & $\theta $ & $a$ & $a$ & $\theta $ \\ 
$b$ & $b$ & $b$ & $\theta $ & $b$ & $b$ & $\theta $ \\ 
$c$ & $c$ & $c$ & $c$ & $\theta $ & $\theta $ & $\theta $ \\ 
$d$ & $d$ & $d$ & $d$ & $d$ & $\theta $ & $\theta $ \\ 
$e$ & $e$ & $e$ & $e$ & $e$ & $e$ & $\theta $%
\end{tabular}%
,
\end{equation*}%
but the block code associated to the BCK-algebra $B$ is \newline
$V_{B}=\{111111,011001,001001,000111,000011,000001\}$ and it is different
from $V_{6}$. Therefore, the above mentioned Theorem needs to be understood
as follows:

i) In some circumstances, we can associate a BCK-algebra to a binary block
code.

ii) The BCK-algebra \thinspace $B$ associated to a block code $V$ generates
the same code $V$ if and only if the order relation $\leq _{s}$, generated
by the skeleton associated to the code $V$, is the same with the order
relation $\preceq $, defined on the obtained BCK-algebra $B$.\medskip

\textbf{Example 6.7.} We consider the code $V_{7}=%
\{11111111,01010101,00111111,00010101,$\newline
$00001111,00000101,00000011,00000001\}=$\newline
$=\{O,X,Y,Z,T,U,V,E\}$, with elements lexicographically ordered. The
associated matrix is

\begin{equation*}
\begin{tabular}{|l|llllllll}
& $O$ & $X$ & $Y$ & $Z$ & $T$ & $U$ & $V$ & $E$ \\ \hline
$O$ & $\mathbf{1}$ & $\mathbf{1}$ & $\mathbf{1}$ & $\mathbf{1}$ & $\mathbf{1}
$ & $\mathbf{1}$ & $\mathbf{1}$ & $\mathbf{1}$ \\ 
$X$ & $0$ & $\mathbf{1}$ & $0$ & $\mathbf{1}$ & $0$ & $\mathbf{1}$ & $0$ & $%
\mathbf{1}$ \\ 
$Y$ & $0$ & $0$ & $\mathbf{1}$ & $\mathbf{1}$ & $\mathbf{1}$ & $\mathbf{1}$
& $\mathbf{1}$ & $\mathbf{1}$ \\ 
$Z$ & $0$ & $0$ & $0$ & $\mathbf{1}$ & $0$ & $\mathbf{1}$ & $0$ & $\mathbf{1}
$ \\ 
$T$ & $0$ & $0$ & $0$ & $0$ & $\mathbf{1}$ & $\mathbf{1}$ & $\mathbf{1}$ & $%
\mathbf{1}$ \\ 
$U$ & $0$ & $0$ & $0$ & $0$ & $0$ & $\mathbf{1}$ & $0$ & $\mathbf{1}$ \\ 
$V$ & $0$ & $0$ & $0$ & $0$ & $0$ & $0$ & $\mathbf{1}$ & $\mathbf{1}$ \\ 
$E$ & $0$ & $0$ & $0$ & $0$ & $0$ & $0$ & $0$ & $\mathbf{1}$%
\end{tabular}%
\end{equation*}

Using the order\ given by the relation $\left( 5.1\right) $, we have that $%
O\preceq X\preceq Z\preceq U\preceq E$, $O\preceq Y\preceq T\preceq V\preceq
E$, $O\preceq Y\preceq Z\preceq U\preceq E,O\preceq Y\preceq T\preceq
U\preceq E$ and the other elements can't be compared. Therefore this order
relation is $\leq _{11}^{8}$and the associated Wajsberg algebra is $\mathcal{%
W}_{11}^{8}$, given by relation $\left( 4.7\right) $.\medskip

\textbf{Example 6.8.} We consider the code $V_{8}=%
\{11111111,01110111,00110011,00010001,$\newline
$00001111,00000111,00000011,00000001\}=$\newline
$=\{O,X,Y,Z,T,U,V,E\}$, with elements lexicographically ordered. Using order
given by the relation $\left( 5.1\right) $, we have that $O\leq X\leq Y\leq
Z\leq E$, $O\leq X\leq Y\leq V\leq E$,\newline
$O\leq X\leq U\leq V\leq E$, $O\leq T\leq U\leq V\leq E$, and the other
elements can't be compared. Therefore this order relation is $\leq _{12}^{8}$
and the associated Wajsberg algebra is $\mathcal{W}_{12}^{8}$, given by
relation $\left( 4.8\right) $.\medskip

\textbf{Example 6.9.} We consider the code $V_{9}=%
\{11111111,01010101,0010001,00010001,$\newline
$00001111,00000101,00000011,00000001\}=$\newline
$=\{O,X,Y,Z,T,U,V,E\}$, with elements lexicographically ordered. Using order
\ given by the relation $\left( 5.1\right) $, we have that $O\leq X\leq
Y\leq Z\leq E$, $O\leq X\leq Y\leq V\leq E$,\newline
$O\leq X\leq U\leq V\leq E$, $O\leq T\leq U\leq V\leq E$, and the other
elements can't be compared. Therefore this order relation is $\leq _{21}^{8}$%
and the associated Wajsberg algebra is $\mathcal{W}_{21}^{8}$, given by the
relation $\left( 4.9\right) $.\medskip

\textbf{Example 6.10. }We consider the code $V_{10}=\{011,101\}$. Using
Theorem 5.3, we have the matrix $B$ 
\begin{equation*}
B=%
\begin{tabular}{|llllllll}
$O$ & $X$ & $Y$ & $Z$ & $T$ & $U$ & $V$ & $E$ \\ 
$1$ & $1$ & $1$ & $1$ & $1$ & $\mathbf{1}$ & $\mathbf{1}$ & $\mathbf{1}$ \\ 
$0$ & $1$ & $0$ & $1$ & $0$ & $\mathbf{1}$ & $\mathbf{0}$ & $\mathbf{1}$ \\ 
$0$ & $0$ & $1$ & $1$ & $0$ & $\mathbf{0}$ & $\mathbf{1}$ & $\mathbf{1}$ \\ 
$0$ & $0$ & $0$ & $1$ & $0$ & $\mathbf{1}$ & $\mathbf{0}$ & $\mathbf{1}$ \\ 
$0$ & $0$ & $0$ & $0$ & $1$ & $\mathbf{1}$ & $\mathbf{1}$ & $\mathbf{1}$ \\ 
$0$ & $0$ & $0$ & $0$ & $0$ & $\mathbf{1}$ & $\mathbf{0}$ & $\mathbf{1}$ \\ 
$0$ & $0$ & $0$ & $0$ & $0$ & $\mathbf{0}$ & $\mathbf{1}$ & $\mathbf{1}$ \\ 
$0$ & $0$ & $0$ & $0$ & $0$ & $\mathbf{0}$ & $\mathbf{0}$ & $\mathbf{1}$%
\end{tabular}%
,
\end{equation*}%
the set $A=\{U,V,E\}$, the Wajsberg algebra $\mathcal{W}_{11}^{8}$ and the
W-function $f:A\rightarrow \mathcal{W}_{11}^{8},f\left( U\right) =U,f\left(
V\right) =V,f\left( E\right) =E$. We have the code $C_{23}=%
\{111,101,011,111,001\}$. The code $V_{10}=\{011,101\}$ is a subcode of the
code $C_{23}$. We remark that the associated Wajsberg is not unique. Indeed,
we can consider the set $A=\{b,c,d\}$, the Wajsberg algebra $\mathcal{W}%
_{12}^{6}$ and the code $C_{23}^{\prime }=\{111,011,101,010,001,000\}$, as
in Example 6.5. The code $V_{10}$ is a subcode of the code $C_{23}^{\prime }$%
.%
\begin{equation*}
\end{equation*}

\textbf{Conclusions}%
\begin{equation*}
\end{equation*}

In this paper, we presented some connections between BCK-commutative bounded
algebras, MV-algebras, Wajsberg algebras and binary block codes. By studying
these connections were identified two types of approaches. First of them is
if these algebras can generate good codes? At a first glance, the answer can
be \textit{no}. Indeed, the codes generates by the BCK-commutative bounded
algebras, MV-algebras and Wajsberg algebras have, in general, the minimum
Hamming distance equal with $1$, that means are not so good codes.
Therefore, we turned our attention to the second approach: if we use the
attached codes, we can find new and interesting properties for these
algebras? The answer is \textit{yes }and that is exactly what we did in this
paper.

We found an algorithm to generate all finite partially ordered Wajsberg
algebras and, from here, all finite partially ordered MV algebras and all
finite partially ordered BCK commutative bounded algebras. In Section 6, we
gave examples of the above mentioned algebras associated to a binary block
codes. We remarked that to a two isomorphic algebras correspond different
codes, but we can find different algebras which can generate the same code.

Even if, for the moment, the answer is \textit{no}, we will not give up the
first approach and we hope that in a future research to achieve important
results in this direction.

\begin{equation*}
\end{equation*}

\textbf{References}%
\begin{equation*}
\end{equation*}

[AAT; 96] Abujabal, H.A.S., Aslam, M., Thaheem, A.B., \textit{A
representation of bounded commutative BCK-algebras}, Internat. J. Math. \&
Math. Sci., 19(4)(1996), 733-736.

[Bu; 06] Bu\c{s}neag, D., \textit{Categories of Algebraic Logic}, Editura
Academiei Rom\^{a}ne, 2006.

[CHA; 58] Chang, C.C.,\textit{\ Algebraic analysis of many-valued logic},
Trans. Amer. Math. Soc. 88(1958), 467-490.

[COM; 00] Cignoli, R. L. O, Ottaviano, I. M. L. D, Mundici, D., \textit{%
Algebraic foundations of many-valued reasoning}, Trends in Logic, Studia
Logica Library, Dordrecht, Kluwer Academic Publishers, 7(2000).

[CT; 96] Cignoli, R., Torell, A., T., \textit{Boolean Products of
MV-Algebras: Hypernormal MV-Algebras}, J Math Anal Appl (199)(1996), 637-653.

[FL; 15] Flaut, C., \textit{BCK-algebras arising from block codes}, Journal
of Intelligent and Fuzzy Systems 28(4)(2015), 1829--1833.

[FL; 17] Flaut, C., \textit{Some Connections Between Binary BlockCodes and
Hilbert Algebras}, in A. Maturo et all, Recent Trends in Social Systems:
Quantitative Theories and Quantitative Models, Springer 2017, p. 249-256.

[FRT; 84] Font, J., M., Rodriguez, A., J., Torrens, A., \textit{Wajsberg
Algebras}, Stochastica, 8(1)(1984), 5-30.

[GA; 90] Gaitan, H., \textit{About quasivarieties of p-algebras and Wajsberg
algebras}, 1990, Retrospective Theses and Dissertations, 9440,
https://lib.dr.iastate.edu/rtd/9440

[HR; 99] H\"{o}hle, U., Rodabaugh, S., E., \textit{Mathematics of Fuzzy
Sets: Logic, Topology and Measure Theory}, Springer Science and Business
Media, LLC, 1999.

[II; 66] Imai, Y., Iseki, K., \textit{On axiom systems of propositional
calculi,} Proc Japan Academic 42(1966), 19--22.

[Io; 08] Iorgulescu, A., Algebras of Logic as BCK Algebras, Editura ASE,
Bucure\c{s}ti, 2008.

[JUN; 11] Jun, Y. B., Song, S. Z., \textit{Codes based on BCK-algebras},
Inform. Sciences., 181(2011), 5102-5109.

[JUN; 03] Jun, Y. B., \textit{Satisfactory filters of BCK-algebras},
Scientiae Mathematicae Japonicae Online, 9(2003), 1--7.

[Me-Ju; 94] Meng, J., Jun, Y. B., \textit{BCK-algebras}, Kyung Moon Sa Co.
Seoul, Korea, 1994.

[Mu; 07] Mundici, D., \textit{MV-algebras-a short tutorial}, Department of
Mathematics \textit{Ulisse Dini}, University of Florence, 2007.

[Pi; 07] Piciu, D., \textit{Algebras of Fuzzy Logic}, Editura Universitaria,
Craiova, 2007.

[WDH; 17] Wang, J., T., Davvaz, B., He, P., F., \textit{On derivations of
MV-algebras}, https://arxiv.org/pdf/1709.04814.pdf%
\begin{equation*}
\end{equation*}

Cristina Flaut,

Faculty of Mathematics and Computer Science,

Ovidius University of Constan\c{t}a, Rom\^{a}nia

cflaut@univ-ovidius.ro

cristina\_flaut@yahoo.com

\bigskip

Radu Vasile,

PhD student at Doctoral School of Mathematics,

Ovidius University of Constan\c{t}a, Rom\^{a}nia

rvasile@gmail.com

\end{document}